\pdfoutput=1
\documentclass[11pt]{article}

\usepackage[margin=1.15in]{geometry}
\usepackage{amsmath,amssymb,amsthm}
\usepackage{graphicx}
\usepackage{subcaption}
\usepackage{booktabs}
\usepackage{microtype}
\usepackage{etoolbox}
\usepackage[hidelinks]{hyperref}

\newtheorem{theorem}{Theorem}[section]
\newtheorem{lemma}[theorem]{Lemma}
\newtheorem{proposition}[theorem]{Proposition}
\newtheorem{corollary}[theorem]{Corollary}
\theoremstyle{definition}
\newtheorem{algorithm}[theorem]{Algorithm}
\theoremstyle{remark}

\newcommand{\sph}{\mathbb{S}^2}
\newcommand{\R}{\mathbb{R}}
\newcommand{\dA}{\,\mathrm{d}A}
\DeclareMathOperator{\tr}{tr}
\DeclareMathOperator{\spann}{span}
\newcommand{\Hperp}{H^1_{\scriptscriptstyle\perp}}

\numberwithin{equation}{section}

\title{\bfseries Solving the Minkowski problem by continuation from the unit sphere: spherical harmonic discretization and error estimates}
\author{Buyang Li \and Jian-Guo Liu \and Rong Tang}
\date{}

\begin{document}

\maketitle

\begin{abstract}
Given a positive function $\rho$ on the unit sphere
$\sph$, the classical Minkowski problem asks for a closed convex
surface in $\R^3$ with Gauss curvature $1/\rho(\nu)$ at the point with
outer normal $\nu$.  In terms of the support function $h$, it is the
Monge--Amp\`ere equation $\det(\nabla^2h+h\sigma)=\rho$ on $\sph$,
uniquely solvable up to translation when $\rho$ has zero first moment.
We solve it by continuation from the sphere: along the linear path
$\rho_t=(1-t)+t\rho$, the solution $h(t)$ evolves in the artificial
time $t$, and its velocity $\partial_th$ solves a linear elliptic
equation in divergence form whose kernel consists of the translations.
In time we use the explicit Euler or trapezoidal method with step
$\tau$; in space, spherical harmonics of degree at most $L$ with a
Lagrange multiplier for the kernel.  For $\rho\in H^{s+2}(\sph)$ with
real $s\ge4$, we prove that all numerical solutions stay strictly
convex, with $H^1$ errors $O(\tau+L^{-(s+1)})$ for Euler and
$O(\tau^3+L^{-(s+1)})$ for the trapezoidal scheme, including data
quadrature, under independent thresholds for $\tau$ and $L$.  The
proof rests on a priori estimates for the perturbed equation satisfied
by each iterate, which keep every step uniformly elliptic, and on the
exactly quadratic determinant, which raises the trapezoidal order
from two to three.  The continuation thus supplies a good initial
guess for Newton's method: from the final state it converges
quadratically to the Galerkin solution, itself within
$O(L^{-(s+1)})$ of the exact solution, whereas started directly from
the sphere its first undamped step left the convex cone in our
ellipsoid test.  Experiments
illustrate the rates, use the multiplier as a compatibility
diagnostic, reconstruct convex polyhedra from facet normals and
areas, and locate umbilics with indices.
\end{abstract}



{\small\setcounter{tocdepth}{1}%
\makeatletter\patchcmd{\l@section}{\addvspace{1.0em \@plus\p@}}{\addvspace{0.5em}}{}{}\makeatother
\tableofcontents}

\section{Introduction}\label{sec:intro}

The classical Minkowski problem asks for a closed convex surface whose
Gauss curvature is prescribed as a function of its outer unit normal.
It arose in Minkowski's work on mixed volumes
\cite{minkowski1903}.  The measure formulation is treated in
\cite{fencheljessen1938} and \cite[Chapter~8]{schneider2014}; modern
extensions include \cite{lutwak1993,hlyz2016}.  If $h$ is the
support function of a smooth strictly convex body in $\R^3$, then
\begin{align*}
W[h]=\nabla^2h+h\sigma,
\qquad
\det W[h]=\rho\quad\text{on }\sph,
\end{align*}
where $\sigma$ is the round metric and $\nabla$ its covariant
derivative. The eigenvalues of $W[h]$ are the principal radii of curvature, and
$\rho$ is the reciprocal Gauss curvature; this is a
Monge--Amp\`ere equation on $\sph$.  Solvability requires
$\rho>0$ and $\int_{\sph}\nu\rho(\nu)\dA(\nu)=0$.
These conditions are sufficient, and the solution is unique up to
translation.  The smooth existence and regularity theory was developed in
\cite{lewy1938,nirenberg1953,pogorelov1978,chengyau1976}; see also
\cite[Chapter~8]{schneider2014}.  We fix the translation by placing the
Steiner point at the origin.

For polyhedral data, Little \cite{little1983} and Lemordant, Pham
Dinh Tao, and Zouaki \cite{lemordanttaozouaki1993} used
finite-dimensional optimization to reconstruct a convex polyhedron from
its extended Gaussian image; related variational principles appear in
\cite{glsy2016}, and Sumbatyan and Troyan \cite{sumbatyantroyan1992}
reconstructed a convex defect from scattered wave data in the ray
approximation.  Lamberg and Kaasalainen \cite{lambergkaasalainen2001}
approximated the support function in finite-dimensional spaces,
including spherical harmonic spaces, and solved the resulting
constrained optimization problem by Newton iterations.  Cheng
\cite{cheng2003} used a level set flow and high-order finite
differences.  Leng, Liu, and Xi \cite{lengliuxi2023} gave a spherical
harmonic reconstruction from prescribed surface area measures, and Hu
and Leng \cite{huleng2025} treated even $L_p$ surface area measures by
finite-dimensional convex optimization; both works prove
convergence of the reconstructed bodies.

We now review the numerical literature on Monge--Amp\`ere
equations.  Geometric discretizations of the planar Dirichlet
problem go back to Oliker and Prussner \cite{olikerprussner1988};
wide-stencil schemes converging to viscosity solutions were developed
in \cite{oberman2008,froeseoberman2011}, the vanishing
moment method in \cite{fengneilan2009}, and quantitative rates for the
Dirichlet problem on flat domains were proved in
\cite{brennergudineilansung2011,nochettozhang2019,nochettontogkaszhang2019rates};
see the survey \cite{neilansalgadozhang2017}.  Brenner
\cite{brenner2026} presents nonlinear least-squares finite element methods
whose degrees of freedom enforce local convexity, with $O(h^2)$ error
estimates for the planar Dirichlet problem, and outlines a least-squares
formulation of the Minkowski problem on the faces of a bounding box as
work in progress.  Dirichlet problems
for graphs of prescribed Gauss curvature on flat domains were solved by
convergent monotone finite difference schemes in
\cite{hamfeldt2018}, and a Dirichlet formulation of the two-dimensional
Minkowski problem by operator splitting and finite elements in
\cite{liuleungqian2024}.  On
the sphere, Hamfeldt and Turnquist
\cite{hamfeldtturnquist2022} proved convergence of
finite difference schemes for optimal transport and reflector
equations.  Kitagawa, M\'erigot, and
Thibert \cite{kmt2019} proved global linear convergence of a damped
Newton iteration for semi-discrete optimal transport.
Kousholt \cite{kousholt2017} reconstructed convex bodies from finitely
many moments of their surface area measures, with a rate in the
translative Hausdorff distance. For the Weyl
isometric embedding problem, whose datum is a metric rather than a
curvature function, a spectral embedding flow was used in
\cite{jasiulekkorzynski2012}, and a recent finite element preprint
proves error estimates with a Lagrange multiplier for the rigid-motion
kernel \cite{gaohulizhang2026}.

We solve the Minkowski problem by continuation from the unit sphere, a
numerical form of the classical continuity method, with the path
parameter $t$ as an artificial time; the problem itself
carries none.
The density $\rho_t=(1-t)+t\rho$, $0\le t\le1$, remains
positive and compatible, and differentiating $\det W[h(t)]=\rho_t$
gives the linear equation
\begin{align*}
 \mathcal{L}_{h(t)}v=\rho-1,\quad v=\partial_th,
 \qquad
 \mathcal{L}_hv=C[h]^{ij}W[v]_{ij}
 =\nabla_i\bigl(C[h]^{ij}\nabla_jv\bigr)+(\tr W[h])\,v,
\end{align*}
where $C[h]=(\tr W[h])\sigma-W[h]$ is the cofactor tensor.
The divergence form holds because $C[h]$ is divergence free, and
$C[h]$ is positive definite when $h$ is strictly convex.  The
continuation thus replaces the nonlinear problem by a sequence of linear
elliptic equations in divergence form, with the symmetric weak form
$b_h(v,q)=\int_{\sph}\bigl(C[h]^{ij}\nabla_jv\,\nabla_iq
-(\tr W[h])\,vq\bigr)\dA$.  The kernel of $\mathcal{L}_h$ is exactly
$\mathcal N=\operatorname{span}\{\nu_1,\nu_2,\nu_3\}$, the
infinitesimal translations, so each stage is one symmetric
saddle-point problem with a Lagrange multiplier $\lambda\in\R^3$:
\begin{align*}
 b_h(v,q)+\sum_{k=1}^3\lambda_k(\nu_k,q)=-(\rho-1,q)
 \quad\text{for all }q,
 \qquad
 (v,\nu_k)=0,\quad k=1,2,3.
\end{align*}
We discretize in space by a spherical harmonic Galerkin method.
The kernel lies in the trial space $V_L$ of spherical polynomials of
degree at most $L$ \cite{atkinsonhan2012}, so the constraints are
imposed without approximation.  Sobolev norms are diagonal in the
spherical harmonic basis \eqref{eq:sobolev}, so approximation and
inverse estimates are explicit.  For trial functions
in $V_L$, the integrands of $b_h$ are spherical polynomials of degree at
most $3L+4$, which the product rule \eqref{eq:QL} integrates exactly
under \eqref{eq:QLexact}; only the data term carries a quadrature
error.

Let $\tau=1/N$, $t_n=n\tau$, and let $\mathcal V_L(h)(\rho-1)$
be the Galerkin velocity at $h\in V_L$.  The Euler scheme sets
$h^{n+1}=h^n+\tau v^n$ with $v^n=\mathcal V_L(h^n)(\rho-1)$; the
explicit trapezoidal scheme sets
$h^{n+1}=h^n+\tfrac{\tau}{2}(k_1+k_2)$ with $k_1=v^n$ and
$k_2=\mathcal V_L(h^n+\tau k_1)(\rho-1)$.  For positive compatible
$\rho\in H^{s+2}(\sph)$ with real $s\ge4$, we prove the following.
{\renewcommand{\labelenumi}{(\roman{enumi})}%
\vspace{-4pt}
\begin{enumerate}
\setlength{\itemsep}{0pt}\setlength{\parsep}{0pt}
\item Error bounds.  The continuation states, without Newton
correction, satisfy
\begin{align*}
 \max_{0\le n\le N}\|h^n-h(t_n)\|_{H^1}\le
 \begin{cases}
  C\bigl(\tau+L^{-(s+1)}\bigr) & \text{(Euler)},\\[2pt]
  C\bigl(\tau^3+L^{-(s+1)}\bigr) & \text{(trapezoidal)},
 \end{cases}
\end{align*}
including the quadrature of the data (Theorems~\ref{thm:main}
and~\ref{thm:trap}); with exact data integrals,
$\rho\in C^{4,\alpha}(\sph)\cap H^s(\sph)$ for some $\alpha\in(0,1)$
suffices.
\item Convexity.  Every continuation state is strictly convex, and the
thresholds $\tau\le\tau_0$ and $L\ge L_0$ are independent of each other.
\item Third order.  The trapezoidal variant is third
order in time because the determinant is exactly quadratic in $h$,
\begin{align*}
 \det W[h+w]=\det W[h]+\mathcal{L}_hw+\det W[w],
\end{align*}
and both schemes use the exact data increment
$\rho_{t_{n+1}}-\rho_{t_n}=\tau(\rho-1)$ of the linear path.
Both features are special to $\sph$ and to this path; see
the discussion before Theorem~\ref{thm:trap}.
\item Compatibility diagnostic.  The multiplier satisfies
$\lambda\cdot\nu=-P_{\mathcal N}(\rho-1)$, with $P_{\mathcal N}$ the
$L^2$-orthogonal projection onto $\mathcal N$ \eqref{eq:PNdef}; for
incompatible data, the continuation flow integrates the projected
density $\rho_t-P_{\mathcal N}\rho_t$ (Proposition~\ref{prop:multiplier}).
\item Newton initialization.  Under $\tau\le cL^{-2}$
(Euler) or $\tau\le cL^{-2/3}$ (trapezoidal), the final continuation
state lies in the $H^1$ ball of radius $\delta_1L^{-2}$ about
$P_Lh(1)$, the degree-$L$ truncation of the exact solution.  Started directly from the sphere instead, Newton's method left the convex cone at its first undamped
step (Section~\ref{subsec:expA}).
\item Newton convergence.  In that ball the Galerkin
problem has a unique solution $h_L^\ast$, strictly convex and within
$O(L^{-(s+1)})$ of $h(1)$, and Newton's method converges to it
quadratically,
$\|h^+-h_L^\ast\|_{H^1}\le CL^2\|h-h_L^\ast\|_{H^1}^2$
(Proposition~\ref{prop:newton}).
\end{enumerate}
\vspace{-4pt}}

The main difficulty is uniform convexity.  The schemes
are explicit, and one step changes $W[h^n]$ by $\tau W[v^n]$, so strict
convexity survives only if $\tau\|W[v^n]\|_{L^\infty}$ is small compared
with the ellipticity constant, and this constant must not deteriorate
over the $N=1/\tau$ steps.  The inverse bound
$\|W[q]\|_{L^\infty}\le CL^2\|q\|_{H^1}$ on $V_L$ only gives
$\|W[v^n]\|_{L^\infty}\le CL^2$, hence the restriction $\tau L^2\le c$,
and a perturbation argument would let the ellipticity constant drift
from step to step.

First, each iterate solves the perturbed Minkowski equation
$\det W[h^n]=\rho_{t_n}+r_n$, and by the quadratic identity above the
defect $r_n$ satisfies exact recursions
(Proposition~\ref{prop:discretedefect}):
\begin{align*}
 r_{n+1}&=r_n+\tau\zeta_n+\tau^2\det W[v^n]
 &&\text{(Euler)},\\
 r_{n+1}&=r_n+\tfrac{\tau}{2}\bigl(\zeta_n+\zeta_n^{\mathrm{st}}\bigr)
 +\tfrac{\tau^2}{4}\det W[k_1-k_2]
 &&\text{(trapezoidal)}.
\end{align*}
Here
$\zeta_n=\bar a_L+(I-P_L)\mathcal{L}_{h^n}(v^n-\hat v^n)$ is the
residual of the Galerkin solve, $\hat v^n$ is the velocity of the
continuous problem at the same state, and $\zeta_n^{\mathrm{st}}$ is the
residual at the stage state.  The residual is the projected quadrature
error plus a term orthogonal to $V_L$, so
$\|\zeta_n\|_{(\Hperp)'}\le CL^{-(s+1)}$ and
$\|\zeta_n\|_{C^{2,\alpha}}\le CL^{3+2\alpha-s}$ (Lemma~\ref{lem:zeta}).

Second, summing the recursion keeps $\|r_n\|_{C^{2,\alpha}}\le\delta$,
so the data $\rho_{t_n}+r_n$ stay bounded below by $\delta$ and bounded in
$C^{2,\alpha}$.  For a strictly convex solution, a maximum principle for
the trace of $W[h^n]$, the positive determinant, and Schauder estimates
then give (Lemma~\ref{lem:apriori})
\begin{align*}
 \kappa_1\sigma\le W[h^n]\le\Lambda_1\sigma,
 \qquad \|h^n\|_{C^{4,\alpha}}\le M_1,
 \qquad \|h^n\|_{H^{s+2}}\le C\bigl(1+\|r_n\|_{H^s}\bigr),
\end{align*}
with constants independent of $n$, $\tau$, and $L$; a
discrete Gronwall argument bounds $\|r_n\|_{H^s}$.  Ellipticity is
therefore not propagated from step to step: once
$\tau\|W[v^n]\|_{L^\infty}\le\kappa_1/2$, the next state satisfies
$W[h^{n+1}]\ge\tfrac{\kappa_1}{2}\sigma$, and the same bounds hold again.

Third, the inverse inequality is applied only to the Galerkin gap
$v^n-P_L\hat v^n$, which is $O(L^{-(s+1)})$ in $H^1$, while $\hat v^n$
is bounded in $C^{2,\alpha}$; hence
$\|W[v^n]\|_{L^\infty}\le C+CL^{1+\alpha-s}\le C$, and
$\tau\le\kappa_1/(2C)$ suffices, independently of $L$.  Finally,
summing the recursion in the dual norm gives
$\|r_n\|_{(\Hperp)'}\le C(\tau+L^{-(s+1)})$, and the stationary
stability estimate
$\|h-h(t)\|_{H^1}\le C\|\det W[h]-\rho_t\|_{(\Hperp)'}$ near the path
(Lemma~\ref{lem:stationary}) turns this into the $H^1$ error bound.

Section~\ref{sec:formulation} collects the geometric
tools, the kernel theorem, the regularity of the path, and the weak
formulation with its multiplier.  Section~\ref{sec:method} defines the
Galerkin scheme and bounds the quadrature error of the data term.
Section~\ref{sec:analysis} proves the stability estimates, the a priori
bounds, the exact defect recursions, the two convergence theorems, and
the convergence of the Newton phase.  Section~\ref{sec:experiments}
reports the experiments: convergence tests on an ellipsoid, a density
without closed-form solution, facet normals and areas of
polyhedra beyond the theory, and umbilic recovery.
The appendix contains the deferred proofs and the umbilic phase
verification.

\section{The Minkowski problem and its continuation formulation}
\label{sec:formulation}

\subsection{Notation and preliminaries}\label{subsec:notation}

Let $\sph=\{x\in\R^3:|x|=1\}$ be the unit sphere with the round metric
$\sigma$ induced from $\R^3$, area element $\dA$, and Levi-Civita covariant
derivative $\nabla$. Latin indices $i,j,l$ refer to tensor components on
$\sph$; indices are raised and lowered with $\sigma$, and repeated indices are
summed. We write $(f,g)=\int_{\sph}fg\dA$ for the $L^2(\sph)$ inner product
and $\|\cdot\|_{H^{p}}$ for the norm of the Sobolev space $H^{p}(\sph)$, defined
for real $p\ge0$ through the spherical harmonic expansion: if
$f=\sum_{l\ge0}f_l$ with $f_l$ in the space $\mathcal{H}_l$ of spherical
harmonics of degree $l$, then
\begin{align}\label{eq:sobolev}
\|f\|_{H^{p}}^2=\sum_{l\ge0}\bigl(1+l(l+1)\bigr)^{p}\|f_l\|_{L^2}^2 .
\end{align}
We write $\|f\|_{W^{1,\infty}}=\|f\|_{L^\infty}+\|\nabla f\|_{L^\infty}$
for Lipschitz functions $f$.
H\"older seminorms on $\sph$ are taken over pairs of
points at geodesic distance at most one, so that
$\|f\|_{C^{k,\alpha}}\le\|f\|_{C^{k,\gamma}}$ for
$0<\alpha\le\gamma$.
Tensor Sobolev and H\"older norms use a fixed finite atlas.
For a symmetric two-tensor field $A$ on $\sph$, $\|A\|_{L^\infty}$ denotes
the essential supremum over $\sph$ of the operator norm of $A$ with respect
to $\sigma$, and inequalities such as $A\ge\kappa\,\sigma$ are understood
in the sense of quadratic forms at every point.
The coordinate functions of $\R^3$ restricted to $\sph$ are denoted
$\nu_1,\nu_2,\nu_3$, so that $\nu=(\nu_1,\nu_2,\nu_3)$ is both the generic
point of $\sph$ and the outer unit normal of the sphere at that point. They
span $\mathcal{H}_1$ and satisfy $\nabla_i\nabla_j\nu_k=-\nu_k\,\sigma_{ij}$.
When spherical coordinates are used, we write
$\nu(\theta,\varphi)=(\sin\theta\cos\varphi,\sin\theta\sin\varphi,
\cos\theta)$, where $0\le\theta\le\pi$ is the polar angle and
$0\le\varphi<2\pi$ is the azimuthal angle; in these
coordinates $\sigma=\mathrm{d}\theta^2+\sin^2\theta\,\mathrm{d}\varphi^2$
and $\dA=\sin\theta\,\mathrm{d}\theta\,\mathrm{d}\varphi$.
We set
\begin{align}\label{eq:kernelspace}
\mathcal{N}=\spann\{\nu_1,\nu_2,\nu_3\}\subset C^\infty(\sph),
\qquad
(\nu_j,\nu_k)=\tfrac{4\pi}{3}\,\delta_{jk},
\end{align}
where $\delta_{jk}$ is the Kronecker symbol, and write
\begin{align}\label{eq:PNdef}
P_{\mathcal{N}}f=\frac{3}{4\pi}\sum_{k=1}^3(f,\nu_k)\,\nu_k
\end{align}
for the $L^2$-orthogonal projection onto $\mathcal{N}$.
We write $f\perp\mathcal N$ when $(f,\nu_k)=0$ for $k=1,2,3$.
Let $\Hperp=\{v\in H^1(\sph):(v,\nu_k)=0,\;k=1,2,3\}$ denote the
$L^2$-orthogonal complement of the kernel space $\mathcal{N}$ in
$H^1(\sph)$.  For a bounded linear functional $F$ on $\Hperp$ we
write
$\|F\|_{(\Hperp)'}=\sup\{|F(q)|/\|q\|_{H^1}:q\in \Hperp\setminus\{0\}\}$,
and analogously with the spectral subspace $V_L^\perp$ of
Section~\ref{subsec:spectral} in place of $\Hperp$.

For a convex body $K\subset\R^3$ containing the origin in its interior, the
support function is $h(\nu)=\max_{y\in K}y\cdot\nu$. If $\partial K$ is
smooth and strictly convex, the point of $\partial K$ with outer normal $\nu$
is
\begin{align}\label{eq:embedding}
r(\nu)=\nabla h(\nu)+h(\nu)\,\nu ,
\end{align}
where $\nabla h$ is the gradient with respect to the round
metric $\sigma$, regarded as a vector of $\R^3$ orthogonal to $\nu$.
In the spherical coordinates $\nu(\theta,\varphi)$ above, for
$0<\theta<\pi$,
\begin{align}\label{eq:gradframe}
\nabla h=\partial_\theta h\,e_\theta
+\frac{\partial_\varphi h}{\sin\theta}\,e_\varphi,
\qquad
e_\theta=\partial_\theta\nu,
\quad
e_\varphi=\frac{\partial_\varphi\nu}{\sin\theta}.
\end{align}
Equivalently, $r(\nu)=\nabla_{\R^3}H(\nu)$ is the Euclidean gradient of
the $1$-homogeneous extension $H(y)=|y|\,h(y/|y|)$
\cite[Section~1.7]{schneider2014}.
Define, for $h\in H^2(\sph)$, the symmetric two-tensor
\begin{align}\label{eq:Wdef}
W[h]_{ij}=\nabla_i\nabla_jh+h\,\sigma_{ij}.
\end{align}
For the support function of a smooth strictly convex body, $W[h]$ is positive
definite and its eigenvalues with respect to $\sigma$ are the principal radii
of curvature at the point $r(\nu)$; consequently
\begin{align*}
\det W[h](\nu)=\frac{1}{\mathcal{K}(r(\nu))},
\end{align*}
where $\mathcal{K}$ is the Gauss curvature and the determinant is
computed with respect to $\sigma$, equivalently in any local
$\sigma$-orthonormal frame \cite[Section~2.5]{schneider2014}.
The map $h\mapsto W[h]$ is linear, and we also write
$W[h]=\nabla^2h+h\sigma$, with $\nabla^2h$ the covariant Hessian.
Furthermore $W[\nu_k]=0$ for $k=1,2,3$: translating $K$ by
$a\in\R^3$ changes $h$ into $h+a\cdot\nu$ without changing curvatures.

The Minkowski problem in its smooth form reads: given $\rho\in C(\sph)$ with
$\rho>0$, find $h$ with $W[h]>0$ such that
\begin{align}\label{eq:minkowski}
\det W[h]=\rho \quad\text{on }\sph .
\end{align}
A solution determines the body up to translation.  We place its
Steiner point
\begin{align}\label{eq:steiner}
\mathrm{st}(h)=\frac{3}{4\pi}\int_{\sph}h(\nu)\,\nu\dA(\nu)\in\R^3 ,
\end{align}
at the origin.  The Steiner point lies in the interior of the body
\cite[Sections~1.7 and~5.4]{schneider2014}; this normalization is equivalent to
$(h,\nu_k)=0$ for $k=1,2,3$.
The following classical
theorem collects the facts we use; see
\cite{lewy1938,nirenberg1953,pogorelov1978,chengyau1976} and
\cite[Chapter~8]{schneider2014}.

\begin{theorem}\label{thm:classical}
Let $\rho\in C^{k,\alpha}(\sph)$ with $k\ge2$, $0<\alpha<1$, $\rho>0$, and
\begin{align}\label{eq:closure}
\int_{\sph}\nu\,\rho(\nu)\dA(\nu)=0 .
\end{align}
Then \eqref{eq:minkowski} has a solution $h\in C^{k+2,\alpha}(\sph)$ with
$W[h]>0$, unique up to addition of an element of $\mathcal{N}$, and hence a
unique solution with $\mathrm{st}(h)=0$; this solution satisfies $h>0$. If
$\rho\in C^\infty(\sph)$ then $h\in C^\infty(\sph)$.
\end{theorem}

We refer to \eqref{eq:closure} as the compatibility condition. It is also
necessary for solvability; this is a special case of the identity in
Lemma~\ref{lem:closureid} below, which holds without any convexity assumption.

\subsection{Cofactor, mixed discriminant, and the trilinear form}
\label{subsec:algebra}

For a symmetric two-tensor $A$ on $\sph$ define the cofactor tensor
\begin{align}\label{eq:cofactor}
C(A)^{ij}=(\tr A)\,\sigma^{ij}-A^{ij},
\end{align}
with $\tr A=\sigma^{ij}A_{ij}$. In an orthonormal frame, $C(A)$ is the usual
cofactor matrix of the two-by-two matrix $A$, so that
\begin{align}\label{eq:cofactorids}
\begin{split}
C(A)^{ij}A_{ij}&=2\det A,\\
\det(A+B)&=\det A+C(A)^{ij}B_{ij}+\det B
\end{split}
\end{align}
for all symmetric $A,B$; the second identity holds because the determinant of
a two-by-two matrix is a quadratic form. The symmetric bilinear form
associated with it is the mixed discriminant
\begin{align}\label{eq:mixeddisc}
D(A,B)=\tfrac12\,C(A)^{ij}B_{ij}
=\tfrac12\bigl(\tr A\,\tr B-A^{ij}B_{ij}\bigr)
=D(B,A),
\qquad D(A,A)=\det A .
\end{align}
For $h\in H^2(\sph)$ we abbreviate $C[h]=C(W[h])$; explicitly,
$C[h]^{ij}=(\Delta h+h)\sigma^{ij}-\nabla^i\nabla^jh$, where
$\Delta=\sigma^{ij}\nabla_i\nabla_j$ is the Laplace--Beltrami operator. If
$W[h]>0$ then $C[h]>0$, since the eigenvalues of $C(A)$ are those of $A$ with
the roles interchanged.

\begin{lemma}[divergence-free cofactor]\label{lem:divfree}
For every $h\in C^3(\sph)$, the tensor $W[h]$ is a Codazzi tensor, that is,
$\nabla_iW[h]_{jk}=\nabla_jW[h]_{ik}$, and consequently
\begin{align}\label{eq:divfree}
\nabla_iC[h]^{ij}=0 .
\end{align}
\end{lemma}

The proof is given in Appendix~\ref{app:two}.
For $h\in H^2(\sph)$, define the linear operator on $v\in H^2(\sph)$ by
\begin{align}\label{eq:Lop}
\mathcal{L}_h v
= C[h]^{ij}\bigl(\nabla_i\nabla_jv+v\,\sigma_{ij}\bigr)
= C[h]^{ij}\,W[v]_{ij}
= 2\,D(W[h],W[v]) ;
\end{align}
by \eqref{eq:cofactorids}, $\mathcal{L}_h$ is the derivative of the map
$h\mapsto\det W[h]$ at $h$, that is, the linearization of the Minkowski
operator. Lemma~\ref{lem:divfree} puts $\mathcal{L}_h$ in divergence
form: expanding
\begin{align*}
\nabla_i\bigl(C[h]^{ij}\nabla_jv\bigr)
=C[h]^{ij}\nabla_i\nabla_jv
+\bigl(\nabla_iC[h]^{ij}\bigr)\nabla_jv
=C[h]^{ij}\nabla_i\nabla_jv
\end{align*}
by \eqref{eq:divfree}, and noting that
$C[h]^{ij}\sigma_{ij}=\tr C[h]=\tr W[h]$ by \eqref{eq:cofactor}, the
definition \eqref{eq:Lop} becomes
\begin{align*}
\mathcal{L}_hv=\nabla_i\bigl(C[h]^{ij}\nabla_jv\bigr)+(\tr W[h])\,v .
\end{align*}
One integration by parts on the closed surface $\sph$ therefore yields,
for $h\in C^3(\sph)$ and $v,q\in H^1(\sph)$ in the weak sense,
\begin{align}\label{eq:greens}
(\mathcal{L}_hv,q) = -b_h(v,q),
\qquad
b_h(v,q):=\int_{\sph}\Bigl(C[h]^{ij}\nabla_jv\,\nabla_iq
-(\tr W[h])\,vq\Bigr)\dA .
\end{align}
Since $C[h]^{ij}$ is a symmetric tensor, the form $b_h$ is symmetric,
$b_h(v,q)=b_h(q,v)$; in this sense $\mathcal{L}_h$ is formally
self-adjoint. Finally, taking $v=h$ in \eqref{eq:Lop} and using
$D(W[h],W[h])=\det W[h]$ from \eqref{eq:mixeddisc} gives Euler's
identity for the quadratic map $h\mapsto\det W[h]$,
\begin{align}\label{eq:Lhh}
\mathcal{L}_hh=2\det W[h].
\end{align}
The trilinear form
\begin{align}\label{eq:trilinear}
T(f,g,w)=\int_{\sph}f\,C[g]^{ij}\,W[w]_{ij}\dA
=2\int_{\sph}f\,D(W[g],W[w])\dA
\end{align}
is symmetric in all three arguments for $f,g,w\in C^2(\sph)$. The
verification is given in Appendix~\ref{app:two}.
For support functions, $\tfrac16T(h_1,h_2,h_3)$ is the mixed volume
$V(K_1,K_2,K_3)$ of the corresponding bodies
\cite[Section~5.1]{schneider2014}, and the symmetry of $T$ is the classical
symmetry of mixed volumes.

\begin{lemma}[compatibility identity]\label{lem:closureid}
For every $h\in C^3(\sph)$,
\begin{align*}
\int_{\sph}\nu_k\,\det W[h]\dA=0,
\qquad k=1,2,3 .
\end{align*}
In particular the compatibility condition \eqref{eq:closure} is necessary
for the solvability of \eqref{eq:minkowski}.
\end{lemma}

\begin{proof}
By \eqref{eq:Lhh}, \eqref{eq:greens}, and $W[\nu_k]=0$,
\begin{align*}
2\int_{\sph}\nu_k\det W[h]\dA
=(\mathcal{L}_hh,\nu_k)
=-b_h(h,\nu_k)
=(\mathcal{L}_h\nu_k,h)
=0 . \tag*{\qedhere}
\end{align*}
\end{proof}

The next lemma is the pointwise Cauchy--Schwarz inequality for mixed
discriminants; it is the two-dimensional case of an inequality of
Alexandrov \cite{alexandrov1938}; see also
\cite[Section~5.5]{schneider2014}.

\begin{lemma}\label{lem:pointwise}
Let $A,B$ be symmetric two-by-two matrices with $A>0$. Then
$D(A,B)^2\ge\det A\,\det B$. If moreover $D(A,B)=0$, then $\det B\le0$, and
$\det B=0$ if and only if $B=0$.
\end{lemma}

The proof is given in Appendix~\ref{app:two}.

\begin{theorem}[kernel of the linearized operator]\label{thm:kernel}
Let $h\in C^3(\sph)$ satisfy $W[h]>0$ and $h>0$. Then
\begin{align*}
\{v\in H^2(\sph):\mathcal{L}_hv=0\}=\mathcal{N}.
\end{align*}
\end{theorem}

The proof is given in Appendix~\ref{app:two}.
Theorem~\ref{thm:kernel} is the infinitesimal form of
Minkowski's uniqueness theorem; compare
\cite[Theorem~8.1.1]{schneider2014}.

\subsection{The continuation path}\label{subsec:path}

Fix a target density $\rho$ satisfying the hypotheses of
Theorem~\ref{thm:classical} and define
\begin{align}\label{eq:path}
\rho_t=(1-t)+t\,\rho ,\qquad t\in[0,1].
\end{align}

We call a density admissible if it is positive and satisfies the
compatibility condition \eqref{eq:closure}, and a data path admissible if
every density along it is.

\begin{lemma}[admissibility of the linear path]\label{lem:path}
For every $t\in[0,1]$, the density $\rho_t$ is positive, satisfies the
compatibility condition $\int_{\sph}\nu\,\rho_t\dA=0$, and
$\rho_0\equiv1$, $\rho_1=\rho$. Moreover
$\partial_t\rho_t=\rho-1$ is independent of $t$, and
$\rho_{t+\tau}=\rho_t+\tau(\rho-1)$ exactly whenever
$t,t+\tau\in[0,1]$.
\end{lemma}

\begin{proof}
Positivity: $\rho_t\ge\min\{1,\min\rho\}>0$ pointwise, as a convex
combination. Compatibility: $\int\nu\dA=0$ and \eqref{eq:closure} give
$\int\nu\rho_t\dA=(1-t)\int\nu\dA+t\int\nu\rho\dA=0$. The remaining
statements are immediate from \eqref{eq:path}.
\end{proof}

By Theorem~\ref{thm:classical}, for each $t\in[0,1]$ there is a unique
$h(t)$ with
\begin{align}\label{eq:pathproblem}
\det W[h(t)]=\rho_t,
\qquad W[h(t)]>0,
\qquad (h(t),\nu_k)=0,\; k=1,2,3,
\end{align}
and $h(0)=1$ is the support function of the unit ball.
Formally differentiating
\eqref{eq:pathproblem} in $t$ and using \eqref{eq:cofactorids} gives the
velocity equation, justified by
Proposition~\ref{prop:pathreg} below: with $v(t)=\partial_th(t)$,
\begin{align}\label{eq:velocity}
\mathcal{L}_{h(t)}v(t)=\rho-1,
\qquad (v(t),\nu_k)=0,\; k=1,2,3 .
\end{align}

\begin{proposition}[regularity of the path]\label{prop:pathreg}
Let $\rho\in C^{k,\alpha}(\sph)$, $k\ge2$, satisfy $\rho>0$ and
\eqref{eq:closure}. Then the family $t\mapsto h(t)$ defined by
\eqref{eq:pathproblem} belongs to
$C^\infty\bigl([0,1];C^{k+2,\alpha}(\sph)\bigr)$, and there are constants
$\kappa_0,\Lambda_0,h_{\min}>0$, depending only on $\rho$, such that
\begin{align}\label{eq:uniformbounds}
\kappa_0\,\sigma\le W[h(t)]\le\Lambda_0\,\sigma,
\qquad h(t)\ge h_{\min},
\qquad t\in[0,1].
\end{align}
\end{proposition}

The proof is given in Appendix~\ref{app:two}.

\subsection{Weak formulation and the multiplier}\label{subsec:weak}

Testing
\eqref{eq:velocity} with $q\in H^1(\sph)$ and using \eqref{eq:greens}, we
arrive at the saddle-point problem: find $v(t)\in H^1(\sph)$ and
$\lambda(t)\in\R^3$ such that
\begin{align}\label{eq:weak}
\begin{aligned}
b_{h(t)}(v(t),q)+\sum_{k=1}^3\lambda_k(t)\,(\nu_k,q)
&=-(\rho-1,\,q)
&&\text{for all }q\in H^1(\sph),\\
(v(t),\nu_k)&=0, && k=1,2,3 .
\end{aligned}
\end{align}
Here $\lambda=(\lambda_1,\lambda_2,\lambda_3)$ is the
Lagrange multiplier of the three constraints, $\lambda_k$ belonging to
$(v,\nu_k)=0$, and we write
$\lambda\cdot\nu=\sum_{k=1}^3\lambda_k\nu_k\in\mathcal N$.
The multiplier has a closed form, and it identifies which equation a
continuation flow integrates when the data path violates the
compatibility condition.

\begin{proposition}[the multiplier measures incompatibility]
\label{prop:multiplier}
Let $h\in C^3(\sph)$ with $W[h]>0$ and $h>0$, and let $f\in L^2(\sph)$.
Consider the problem: find $(v,\lambda)\in H^1(\sph)\times\R^3$ such
that
\begin{align*}
b_h(v,q)+\sum_{k=1}^3\lambda_k\,(\nu_k,q)=-(f,q)
\quad\text{for all }q\in H^1(\sph),
\qquad
(v,\nu_k)=0,\; k=1,2,3.
\end{align*}
Problem~\eqref{eq:weak} is the special case
$h=h(t)$, $f=\rho-1$.

(i) Well-posedness and the multiplier formula. The problem has
a unique solution; the multiplier is explicit,
\begin{align}\label{eq:lambdaformula}
\lambda_k=-\frac{3}{4\pi}\,(f,\nu_k),
\qquad k=1,2,3,
\qquad\text{equivalently}\quad
\lambda\cdot\nu=-P_{\mathcal{N}}f;
\end{align}
and $v$ satisfies $\mathcal{L}_hv=f+\lambda\cdot\nu$ in the weak
sense.

(ii) The flow integrates the projected data path. Let
$t\mapsto\bar\rho_t$ be continuously differentiable into $L^2(\sph)$
with $\bar\rho_0=1$, and let $t\mapsto\bar h(t)$ be continuously
differentiable into $H^2(\sph)$ with $\bar h(t)\in C^3(\sph)$,
$W[\bar h(t)]>0$ and $\bar h(t)>0$ for every $t$, and
$\det W[\bar h(0)]=1$. If for every $t$ the pair
$(\partial_t\bar h(t),\lambda(t))$ solves the problem above with
$h=\bar h(t)$ and $f=\partial_t\bar\rho_t$, then
\begin{align}\label{eq:wrongtarget}
\det W[\bar h(t)]
=\bar\rho_t-P_{\mathcal{N}}\bar\rho_t ,
\qquad t\in[0,1].
\end{align}
In particular, for the linear path \eqref{eq:path} with compatible
target, $\lambda(t)\equiv0$.
\end{proposition}

The proof is given in Appendix~\ref{app:two}.
Define, for $0<\kappa\le\Lambda$, the coefficient tube
\begin{align*}
\mathcal{W}(\kappa,\Lambda)
=\bigl\{h\in C^3(\sph):
\kappa\,\sigma\le W[h]\le\Lambda\,\sigma\bigr\}.
\end{align*}

\begin{lemma}[G\aa rding inequality]\label{lem:garding}
For $h\in\mathcal{W}(\kappa,\Lambda)$ and all $v\in H^1(\sph)$,
\begin{align*}
b_h(v,v)\ \ge\ \kappa\,\|v\|_{H^1}^2-(\kappa+2\Lambda)\,\|v\|_{L^2}^2 .
\end{align*}
\end{lemma}

\begin{proof}
The eigenvalues of the cofactor are those of $W[h]$ interchanged, so
membership in the tube gives the two pointwise bounds
\begin{align*}
C[h]\ \ge\ \kappa\,\sigma,
\qquad
\tr W[h]\ \le\ 2\Lambda .
\end{align*}
Inserting them into the two terms of $b_h$,
\begin{align*}
b_h(v,v)
=\int_{\sph}\Bigl(C[h]^{ij}\nabla_jv\,\nabla_iv-(\tr W[h])\,v^2\Bigr)\dA
\ \ge\ \kappa\,\|\nabla v\|_{L^2}^2-2\Lambda\,\|v\|_{L^2}^2 ,
\end{align*}
and $\|\nabla v\|_{L^2}^2=\|v\|_{H^1}^2-\|v\|_{L^2}^2$ finishes the
proof.
\end{proof}

\begin{proposition}[inf-sup stability on the constrained space]
\label{prop:infsup}
Let $h\in C^3(\sph)$ with $W[h]>0$ and $h>0$. Then there is
$\beta=\beta(h)>0$ such that
\begin{align}\label{eq:infsup}
\inf_{v\in \Hperp\setminus\{0\}}\;
\sup_{q\in \Hperp\setminus\{0\}}\;
\frac{b_h(v,q)}{\|v\|_{H^1}\|q\|_{H^1}}\;\ge\;\beta .
\end{align}
By the Banach--Ne\v{c}as--Babu\v{s}ka theorem and the symmetry of
$b_h$, \eqref{eq:infsup} says exactly that for every bounded linear
functional $F$ on $\Hperp$ the problem
\begin{align*}
v\in \Hperp:
\qquad
b_h(v,q)=F(q)
\quad\text{for all }q\in \Hperp
\end{align*}
has a unique solution, with
$\|v\|_{H^1}\le\beta^{-1}\|F\|_{(\Hperp)'}$.
Moreover:

(i) Uniformity along the path. The constant $\beta$ can be chosen uniformly
over the compact family $\{h(t):t\in[0,1]\}$ of
Proposition~\ref{prop:pathreg}; from now on, $\beta$ denotes
this uniform constant.

(ii) Robustness under coefficient perturbations. Let
$\tilde h\in C^3(\sph)$ satisfy
$\|W[\tilde h]-W[h(t)]\|_{L^\infty}\le\varepsilon$
for some $t\in[0,1]$.
Then
\begin{align*}
\inf_{v\in \Hperp\setminus\{0\}}\;
\sup_{q\in \Hperp\setminus\{0\}}\;
\frac{b_{\tilde h}(v,q)}{\|v\|_{H^1}\|q\|_{H^1}}
\;\ge\;\beta-2\varepsilon.
\end{align*}
In particular, the constrained problems for $b_{\tilde h}$ remain well
posed whenever $\varepsilon<\beta/2$.
This form is used in Section~\ref{sec:analysis}, where the discrete
iterates lie near but not on the continuous path.
\end{proposition}

The proof is given in Appendix~\ref{app:two}.

\section{The numerical method}\label{sec:method}

\subsection{Spectral Galerkin discretization}\label{subsec:spectral}

For an integer $L\ge1$ let
\begin{align*}
V_L=\bigoplus_{l=0}^{L}\mathcal{H}_l,
\qquad \dim V_L=(L+1)^2,
\qquad
V_L^\perp=\{q\in V_L:(q,\nu_k)=0,\;k=1,2,3\},
\end{align*}
and let $P_L:L^2(\sph)\to V_L$ be the $L^2$-orthogonal projection, which is
simultaneously the $H^{p}$-orthogonal projection for every $p$ by
\eqref{eq:sobolev}. Since $\mathcal{N}=\mathcal{H}_1\subset V_L$, the
kernel of the continuous problem is represented exactly, and
$V_L^\perp=V_L\cap \Hperp$. Two standard properties follow directly from
\eqref{eq:sobolev}: the approximation property
\begin{align}\label{eq:approx}
\|f-P_Lf\|_{H^{p'}}\le\bigl(1+L(L+1)\bigr)^{(p'-p)/2}\|f\|_{H^{p}},
\qquad 0\le p'\le p,
\end{align}
and the inverse inequality
\begin{align}\label{eq:inverse}
\|q\|_{H^{p}}\le\bigl(1+L(L+1)\bigr)^{(p-p')/2}\|q\|_{H^{p'}},
\qquad q\in V_L,\quad 0\le p'\le p .
\end{align}
We also use the Nikolskii inequality on $\sph$: every
spherical polynomial $q\in\bigoplus_{l\le L+2}\mathcal{H}_l$
satisfies
\begin{align}\label{eq:nikolskii}
\|q\|_{L^\infty}\le\frac{L+3}{\sqrt{4\pi}}\,\|q\|_{L^2},
\end{align}
and the constant is sharp \cite{atkinsonhan2012}.  A homogeneous-extension argument and
\eqref{eq:nikolskii} give a degree count and inverse-type pointwise
bounds on $V_L$.

\begin{lemma}[degree count and pointwise bounds on $V_L$]\label{lem:Winfty}
Let $q\in V_L$.
 In any Cartesian frame, the entries of $W[q]$ are
spherical polynomials of degree at most $L+2$, and those of
$\nabla q+q\,\nu$ of degree at most $L+1$.  Moreover,
\begin{align}\label{eq:Winfty}
\|W[q]\|_{L^\infty}\le C\,L^{2}\,\|q\|_{H^1},
\qquad
\|q\|_{W^{1,\infty}}\le C\,L\,\|q\|_{H^1}.
\end{align}
\end{lemma}

The proof is given in Appendix~\ref{app:three}.

\subsection{The fully discrete continuation scheme}\label{subsec:scheme}

Let $N\in\mathbb{N}$ be the number of time steps, $\tau=1/N$ and
$t_n=n\tau$; the step number $N$ is not to be confused with the kernel
space $\mathcal{N}$ of \eqref{eq:kernelspace}. The discrete velocity
problem is the saddle-point problem \eqref{eq:weak} posed on $V_L$,
with the data integral computed by the quadrature rule $Q_L$ of
\eqref{eq:QL} below:
given a state $h\in V_L$ and a source $f\in C(\sph)$, find
$(v,\lambda)\in V_L\times\R^3$ such that
\begin{align}\label{eq:discretevelocity}
\begin{aligned}
b_h(v,q)+\sum_{k=1}^3\lambda_k\,(\nu_k,q)&=-Q_L(fq)
&&\text{for all }q\in V_L,\\
(v,\nu_k)&=0,&& k=1,2,3 .
\end{aligned}
\end{align}
We write $v=\mathcal{V}_L(h)f$ for its first component whenever
\eqref{eq:discretevelocity} is uniquely solvable;
Lemma~\ref{lem:discreteinfsup} guarantees this for $L\ge L_0$ and
$h$ near the exact path.

\begin{algorithm}[continuation with Newton correction]\label{alg:scheme}
Input: $\rho$, degree $L$, step number $N$, correction steps $M\ge0$. Set
$h^0=1$ and $f=\rho-1$.
\begin{enumerate}
\item (continuation phase) For $n=0,\dots,N-1$: compute
$v^n=\mathcal{V}_L(h^n)f$ and set either
\begin{align*}
h^{n+1}=h^n+\tau v^n
\end{align*}
in the Euler variant, or
\begin{align*}
k_1=v^n,\qquad
g^n=h^n+\tau k_1,\qquad
k_2=\mathcal{V}_L(g^n)f,\qquad
h^{n+1}=h^n+\tfrac{\tau}{2}(k_1+k_2)
\end{align*}
in the trapezoidal variant, with $g^n$ the stage state.  The stage
quantities $k_1,k_2$ depend on $n$ as well; that index is suppressed.
\item (Newton phase) Set $h^{N,0}=h^N$ and, for $m=1,\dots,M$, compute
$h^{N,m}=h^{N,m-1}
+\mathcal{V}_L(h^{N,m-1})\bigl(\rho-\det W[h^{N,m-1}]\bigr)$.
\end{enumerate}
Output: the final iterate, denoted again by $h^N$ when $M=0$ and by
$h^{N,M}$ otherwise, together with the multipliers $\lambda^n\in\R^3$
of the solves $v^n=\mathcal V_L(h^n)f$.
\end{algorithm}

Every iterate satisfies $(h^n,\nu_k)=0$: the initial value $h^0=1$
satisfies the constraint, and every update lies in $V_L^\perp$ by the
second equation of \eqref{eq:discretevelocity}. By \eqref{eq:steiner}
the Steiner point of every iterate is exactly at the origin.

Both variants of Algorithm~\ref{alg:scheme} rest on the quadratic
expansion of the Minkowski operator: by \eqref{eq:cofactorids} and
the linearity of $h\mapsto W[h]$,
\begin{align}\label{eq:quadratic}
\det W[g+w]=\det W[g]+\mathcal{L}_gw+\det W[w],
\qquad g,w\in H^2(\sph),
\end{align}
with no remainder. By \eqref{eq:Lop} and \eqref{eq:mixeddisc} the
middle term equals $2D(W[g],W[w])$ and is therefore symmetric,
$\mathcal{L}_gw=\mathcal{L}_wg$.

All integrals are computed by the positive tensor-product rule
\begin{align}\label{eq:QL}
Q_L(F)=\frac{2\pi}{n_\varphi}
\sum_{i=1}^{n_x}w_i\sum_{j=0}^{n_\varphi-1}
F\bigl(\sqrt{1-x_i^2}\cos\varphi_j,
\sqrt{1-x_i^2}\sin\varphi_j,x_i\bigr),
\end{align}
with Gauss--Legendre nodes $x_i$ and weights $w_i$ on $[-1,1]$ and
$\varphi_j=2\pi j/n_\varphi$.  Since $\mathrm{d}A=\mathrm{d}x\,\mathrm{d}\varphi$
for $x=\nu_3$, $Q_L$ is exact for spherical polynomials of degree at
most $E_L:=\min\{2n_x-1,\,n_\varphi-1\}$: the azimuthal sum is exact for
$e^{\mathrm{i}m\varphi}$ with $|m|<n_\varphi$ and leaves the zonal part,
a polynomial in $x$ of no higher degree, which Gauss--Legendre
integrates exactly.  The form and constraint terms of
\eqref{eq:discretevelocity} have integrands of degree at most $3L+4$ by
Lemma~\ref{lem:Winfty}, so $Q_L$ computes them exactly, as in
\eqref{eq:discretevelocity}, provided that
\begin{align}\label{eq:QLexact}
E_L\ge3L+4,
\qquad\text{that is,}\qquad
n_x\ge\bigl\lceil\tfrac{3L+5}{2}\bigr\rceil,
\quad n_\varphi\ge3L+5 .
\end{align}
Only the data term $Q_L(fq)$ then differs from the exact
integral; for $f=\rho-1$, the next proposition bounds the difference.

\begin{proposition}[quadrature error of the data term]\label{prop:quadratureload}
Let $s\ge4$ and $\rho\in H^{s+2}(\sph)$.  Define $a_L\in V_L$ by
\begin{align}\label{eq:loadrepresenter}
(a_L,q)=Q_L\bigl((\rho-1)q\bigr)-(\rho-1,q),
\qquad q\in V_L .
\end{align}
Then
\begin{align}\label{eq:loadbounds}
\begin{aligned}
\|a_L\|_{L^2}&\le C L^{-(s+1)}\|\rho\|_{H^{s+2}},
&\|a_L\|_{H^p}&\le C L^{p-s-1}\|\rho\|_{H^{s+2}},
\quad 0\le p\le s,\\
\|a_L\|_{(\Hperp)'}&\le C L^{-(s+1)}\|\rho\|_{H^{s+2}},
&\|a_L\|_{C^{2,\alpha}}&\le
C L^{2+2\alpha-s}\|\rho\|_{H^{s+2}}
\end{aligned}
\end{align}
for $0<\alpha<1/4$.  The constants are independent of $L$.
\end{proposition}

The proof is given in Appendix~\ref{app:three}.

\section{Error analysis}\label{sec:analysis}

Throughout this section, assume that $\rho>0$ satisfies the
compatibility condition \eqref{eq:closure} and
\begin{align*}
\rho\in H^{s+2}(\sph),\qquad s\ge4,
\end{align*}
where $s$ may be noninteger.  Fix $0<\alpha_0<1$.
The embedding $H^{s+2}(\sph)\hookrightarrow H^6(\sph)
\hookrightarrow C^{4,\alpha_0}(\sph)$ allows us to apply
Theorem~\ref{thm:classical} and Proposition~\ref{prop:pathreg}
with $(k,\alpha)=(4,\alpha_0)$.
Let $a_L$ be the quadrature-error representer of
Proposition~\ref{prop:quadratureload} and set
\begin{equation}\label{eq:abar}
\bar a_L:=(I-P_{\mathcal N})a_L.
\end{equation}%
If the data integrals are evaluated exactly, then $a_L=0$ and
the weaker assumption
$\rho\in C^{4,\alpha_0}(\sph)\cap H^s(\sph)$ suffices.

Three velocities appear below.  At a state $h\in V_L$, the
computed velocity $v=\mathcal V_L(h)(\rho-1)$ solves the discrete
problem \eqref{eq:discretevelocity}.  At the same state,
$\hat v\in\Hperp$ denotes the solution of the linearized problem
\begin{equation}\label{eq:vhat}
\mathcal L_h\hat v=\rho-1,
\end{equation}%
with neither Galerkin truncation nor quadrature, whenever this problem
is well posed; the difference $v-\hat v$ is thus the error of the space
discretization at the fixed state $h$.  The path velocity
$v(t)=\partial_th(t)$ of \eqref{eq:velocity} solves the same equation
at the exact state $h(t)$ instead of $h$.  Besides the convexity
bounds \eqref{eq:uniformbounds}, the exact path satisfies
\begin{align}\label{eq:pathregularity}
\sup_{t\in[0,1]}\bigl(\|h(t)\|_{H^{s+2}}+\|v(t)\|_{H^{s+2}}\bigr)<\infty ,
\end{align}
which follows from Lemma~\ref{lem:apriori} with $\psi=\rho_t$ for $h(t)$
and from Lemma~\ref{lem:lineartame} for $v(t)$; see Appendix~\ref{app:four}.  Unless
stated otherwise, $C$ denotes a generic constant depending only on
$\rho$, $s$, $\alpha_0$, the exponent $\alpha$ fixed in
\eqref{eq:alphachoice}, and previously fixed constants, but not on $L$,
$\tau$, or $n$.

\subsection{Preliminary estimates}

\begin{lemma}[coefficient perturbation]\label{lem:perturbation}
Let $g,e,\varphi\in C^2(\sph)$ and $q\in H^1(\sph)$. Then
\begin{align}\label{eq:bsym}
b_{g+e}(\varphi,q)-b_g(\varphi,q)
=b_e(\varphi,q)=b_\varphi(e,q),
\end{align}
and hence
\begin{align}\label{eq:bbound}
|b_{g+e}(\varphi,q)-b_g(\varphi,q)|
\le C\|\varphi\|_{C^2}\|e\|_{H^1}\|q\|_{H^1}.
\end{align}
\end{lemma}

The proof is given in Appendix~\ref{app:four}.

\begin{lemma}[higher-order elliptic estimate]\label{lem:lineartame}
Let $s\ge4$, $0<\alpha<1$, and let $h\in H^{s+2}(\sph)$ satisfy
\begin{align*}
\kappa\sigma\le W[h]\le\Lambda\sigma,
\qquad \|h\|_{C^{4,\alpha}}\le M_0 .
\end{align*}
Assume that the inf-sup bound \eqref{eq:infsup} holds at $h$
with $\beta/2$ in place of $\beta$.
Then, for every $f\in H^s(\sph)$ with $f\perp\mathcal N$, the problem
\begin{align*}
\mathcal L_h w=f,\qquad w\perp\mathcal N,
\end{align*}
has a unique solution $w\in H^{s+2}(\sph)$ satisfying
\begin{align}\label{eq:lineartame}
\begin{aligned}
\|w\|_{C^{2,\alpha}}
&\le C\|f\|_{C^\alpha},\\
\|w\|_{H^{s+2}}
&\le C\left(
\|f\|_{H^s}+
\bigl(1+\|h\|_{H^{s+2}}\bigr)\|f\|_{C^\alpha}
\right),
\end{aligned}
\end{align}
where $C$ depends only on $s,\alpha,\kappa,\Lambda,M_0$, and $\beta$.
\end{lemma}

The proof is given in Appendix~\ref{app:four}.

\begin{lemma}[uniform discrete stability]\label{lem:discreteinfsup}
There exist an integer $L_0\ge1$ and constants
$\varepsilon_0,\beta_0>0$ with the following property.
Let $L\ge L_0$ and $h\in C^3(\sph)$ satisfy
\(\|W[h]-W[h(t)]\|_{L^\infty}\le\varepsilon_0\)
for some $t\in[0,1]$.
Then
\begin{align}\label{eq:discreteinfsup}
\inf_{v\in V_L^\perp\setminus\{0\}}\;
\sup_{q\in V_L^\perp\setminus\{0\}}\;
\frac{b_{h}(v,q)}{\|v\|_{H^1}\|q\|_{H^1}}
\ge\beta_0 .
\end{align}

Consequently, for every $f\in C(\sph)$, the discrete system
\eqref{eq:discretevelocity} at $h$ has a unique solution
$(v,\lambda)$ with
\(\|v\|_{H^1}+|\lambda|\le C\,Q_L(f^2)^{1/2}\).
\end{lemma}

The proof is given in Appendix~\ref{app:four}.

\begin{corollary}[discrete multiplier]\label{cor:lambda}
Let $a_L\in V_L$ be the quadrature-error representer from
\eqref{eq:loadrepresenter}. For fixed $L$, every well-posed
solve of \eqref{eq:discretevelocity} in Algorithm~\ref{alg:scheme}, in
either phase, has the same multiplier,
\[
\lambda_k
=-\frac{3}{4\pi}Q_L(\rho\nu_k)
=-\frac{3}{4\pi}\bigl[(\rho-1,\nu_k)+(a_L,\nu_k)\bigr].
\]
For compatible data,
\[
\lambda\cdot\nu=-P_{\mathcal N}a_L,
\qquad
|\lambda|\le C L^{-(s+1)}\|\rho\|_{H^{s+2}}.
\]
In particular, $\lambda=0$ for compatible data when the data integrals are exact.
\end{corollary}

\begin{proof}
Exactly as for \eqref{eq:lambdaformula}, testing
\eqref{eq:discretevelocity} with $q=\nu_k$ gives
$\lambda_k=-\frac{3}{4\pi}Q_L(f\nu_k)$.
In the continuation phase $f=\rho-1$; in the Newton phase
$f=\rho-\det W[g]$ with $g\in V_L$. Moreover,
\[
Q_L(\nu_k)=0,
\qquad
Q_L(\det W[g]\nu_k)
=(\det W[g],\nu_k)=0,
\]
where quadrature is exact because
$\deg(\det W[g]\nu_k)\le2L+5$ by Lemma~\ref{lem:Winfty},
and the integral vanishes by Lemma~\ref{lem:closureid}.
Hence both the continuation and Newton phases give the stated formula.
For compatible data, $(\rho-1,\nu_k)=0$, so
$\lambda_k=-\frac{3}{4\pi}(a_L,\nu_k)$, that is,
$\lambda\cdot\nu=-P_{\mathcal N}a_L$; since
$\|\nu_k\|_{L^2}^2=4\pi/3$, the Cauchy--Schwarz inequality and the
$L^2$ bound in \eqref{eq:loadbounds} give
$|\lambda|\le CL^{-(s+1)}\|\rho\|_{H^{s+2}}$.
\end{proof}

\begin{lemma}[stationary stability]\label{lem:stationary}
There exist $\varepsilon_1,C>0$, depending only on the exact path,
such that, for every $t\in[0,1]$ and $h\in C^3(\sph)$ with
$h\perp\mathcal N$ and
\[
\|W[h]-W[h(t)]\|_{L^\infty}\le\varepsilon_1,
\]
we have
\begin{equation*}
\|h-h(t)\|_{H^1}
\le C\|\det W[h]-\rho_t\|_{(\Hperp)'}.
\end{equation*}
\end{lemma}

\begin{proof}
Set $e=h-h(t)\in \Hperp$ and
$h_{1/2}=\frac12(h+h(t))$. By \eqref{eq:quadratic} and
\eqref{eq:Lhh},
\[
r:=\det W[h]-\rho_t
=\mathcal L_{h(t)}e+\det W[e]
 =\mathcal L_{h_{1/2}}e.
\]
Choose $0<\varepsilon_1\le\beta/2$.
Since $W[h_{1/2}]-W[h(t)]=\tfrac12W[e]$, the hypothesis gives
\[
\|W[h_{1/2}]-W[h(t)]\|_{L^\infty}\le\tfrac12\varepsilon_1,
\]
so
Proposition~\ref{prop:infsup}(ii), applied with
$\varepsilon=\varepsilon_1/2$, gives an inf-sup lower bound
$\beta-\varepsilon_1\ge\beta/2$ at $h_{1/2}$. Using
$r=\mathcal L_{h_{1/2}}e$ and \eqref{eq:greens}, we obtain
\[
\frac{\beta}{2}\|e\|_{H^1}
\le
\sup_{q\in \Hperp\setminus\{0\}}
\frac{|b_{h_{1/2}}(e,q)|}{\|q\|_{H^1}}
=
\|r\|_{(\Hperp)'}. \qedhere
\]
\end{proof}

\subsection{A priori and residual bounds}
\label{subsec:apriori}

To close the discrete induction and control spectral errors,
we derive ellipticity and regularity directly from the equation
satisfied by each iterate.
With the equation defect $r_n:=\det W[h^n]-\rho_{t_n}$, the
iterate $h^n$ of Algorithm~\ref{alg:scheme} solves the perturbed
Minkowski equation
\begin{equation}\label{eq:exactperturbed}
 \det W[h^n]=\psi_n:=\rho_{t_n}+r_n.
\end{equation}
For strictly convex iterates, Lemma~\ref{lem:apriori} below obtains all
quantitative bounds from $\psi_n$, with constants independent of
$n,\tau,L$ and without inverse estimates on $h^n$.  Set
\begin{equation}\label{eq:deltalambda}
 \delta=\frac14\min\{1,\min_{\sph}\rho\},
 \qquad
 \gamma=\min\{\alpha_0,\tfrac14\},
\end{equation}
where $\alpha_0$ is fixed at the beginning of
Section~\ref{sec:analysis}, and choose
\begin{equation}\label{eq:Lambdastar}
 \Lambda_*\ge
 1+\delta+
 \sup_{t\in[0,1]}
 \|\rho_t\|_{C^{2,\gamma}}.
\end{equation}
The supremum is finite, since
$\|\rho_t\|_{C^{2,\gamma}}\le1+\|\rho\|_{C^{2,\gamma}}$.
With $\alpha_*$ from Lemma~\ref{lem:apriori} below, fix
\begin{equation}\label{eq:alphachoice}
 0<\alpha<
 \min\{\gamma,\alpha_*(\delta,\Lambda_*)\}.
\end{equation}
Since $\rho_t\ge4\delta$ and
$\|\rho_t\|_{C^{2,\alpha}}\le\|\rho_t\|_{C^{2,\gamma}}$ for $\alpha<\gamma$,
\[
 \|r_n\|_{C^{2,\alpha}}\le\delta
 \quad\Longrightarrow\quad
 \psi_n\ge3\delta,
 \qquad
 \|\psi_n\|_{C^{2,\alpha}}\le\Lambda_*,
\]
so Lemma~\ref{lem:apriori} applies with $(\psi_-,\psi_+)=(\delta,\Lambda_*)$
whenever $W[h^n]>0$, since $h^n\in V_L^\perp$ is smooth and
normalized.

\begin{lemma}[a priori bounds from the equation]
\label{lem:apriori}
For $0<\psi_-\le\psi_+$, there exists
$\alpha_*=\alpha_*(\psi_-,\psi_+)\in(0,1)$
such that the following holds.
Let $0<\alpha<\alpha_*$ and
$h\in C^4(\sph)$ satisfy
\begin{equation}\label{eq:perturbeddata}
 W[h]>0,
 \qquad
 (h,\nu_k)=0\quad(k=1,2,3),
 \qquad
 \det W[h]=\psi,
\end{equation}
with
\begin{equation}\label{eq:databounds}
 \psi\ge\psi_-,
 \qquad
 \|\psi\|_{C^{2,\alpha}}\le\psi_+.
\end{equation}
Then
\begin{equation}\label{eq:corridorA}
 \kappa_1\sigma\le W[h]\le \Lambda_1\sigma,
 \qquad
 h\ge\kappa_1,
 \qquad
 \|h\|_{C^{4,\alpha}}\le M_1,
\end{equation}
where $\kappa_1,\Lambda_1,M_1>0$ depend only on
$\psi_-,\psi_+,\alpha$.
If also $\psi\in H^p(\sph)$ for real $p\in[4,s]$, then
\begin{equation}\label{eq:tame}
 \|h\|_{H^{p+2}}
 \le C(\psi_-,\psi_+,\alpha,p)
       \bigl(1+\|\psi\|_{H^p}\bigr).
\end{equation}
\end{lemma}

The hypothesis $W[h]>0$ is qualitative; no lower curvature
bound enters as input.  Once one time step preserves strict convexity,
the equation and the bounds \eqref{eq:databounds} on $\psi$ restore the
corridor \eqref{eq:corridorA}.

\begin{proof}
\emph{Step 1: ellipticity and low-order bounds.}
Set $A=W[h]$, $H=\tr A$, and
\begin{equation*}
 F(A):=\sqrt{\det A},
 \qquad
 F(W[h])=\sqrt{\psi},
 \qquad
 F^{ij}:=\frac{\partial F}{\partial A_{ij}}
         =\frac{C(A)^{ij}}{2\sqrt{\psi}}.
\end{equation*}
Since $A>0$, the tensor $F^{ij}$ is positive definite, and bounding
$H$ bounds the largest eigenvalue of $A$.  Three identities are
available: by Lemma~\ref{lem:divfree}, $A$ satisfies the
Codazzi identity $\nabla_iA_{jk}=\nabla_jA_{ik}$, so it
obeys the Simons identity $\Delta A=\nabla^2H+2A-H\sigma$ of
\eqref{eq:simons} proved in Appendix~\ref{app:four}; $F$ is positively homogeneous of degree one, so
$F^{ij}A_{ij}=\sqrt{\psi}$; and $C(A)^{ij}\sigma_{ij}=H$ by
\eqref{eq:cofactor}, so $F^{ij}\sigma_{ij}=H/(2\sqrt{\psi})\ge1$ by
$H\ge2\sqrt{\det A}$.  Since $h\in C^4$, the tensor $A$ is twice
differentiable, and the concavity of $F$ gives
$F^{ij}\Delta A_{ij}\ge\Delta F(A)=\Delta\sqrt{\psi}$.
Contracting the Simons identity with $F^{ij}$ and inserting the three
identities and this inequality gives
\[
 F^{ij}\nabla_i\nabla_jH
 =F^{ij}\Delta A_{ij}-2\sqrt{\psi}+HF^{ij}\sigma_{ij}
 \ge\Delta\sqrt{\psi}-2\sqrt{\psi}+H .
\]
At a maximum point of $H$ we have $F^{ij}\nabla_i\nabla_jH\le0$,
because $\nabla^2H\le0$ and $F^{ij}>0$.  Therefore
\begin{equation}\label{eq:uppercurv}
 \|H\|_{L^\infty}
 \le2\|\sqrt{\psi}\|_{L^\infty}+\|\Delta\sqrt{\psi}\|_{L^\infty}
 \le \Lambda_1,
 \qquad
 W[h]\le \Lambda_1\sigma,
\end{equation}
with $\Lambda_1=\Lambda_1(\psi_-,\psi_+)$, since
$\Delta\sqrt{\psi}=\Delta\psi/(2\sqrt{\psi})-|\nabla\psi|^2/(4\psi^{3/2})$
and
$\psi\ge\psi_-$.  The two eigenvalues of $A$ multiply to
$\psi\ge\psi_-$, so
\begin{equation}\label{eq:lowercurv}
 W[h]\ge\frac{\psi_-}{\Lambda_1}\sigma
 =:\kappa_1\sigma.
\end{equation}

The trace bound controls $h$ through $(\Delta+2)h=H$: the
normalization \eqref{eq:perturbeddata} places $h$ in the
orthogonal complement of
$\ker(\Delta+2)=\operatorname{span}\{\nu_1,\nu_2,\nu_3\}$, on which
$\Delta+2$ is invertible with inverse bounded from $L^4(\sph)$ to
$W^{2,4}(\sph)$.  Elliptic regularity and Sobolev embedding
\cite[Chapters~7 and~9]{gilbargtrudinger2001} give
\begin{equation*}
 \|h\|_{C^1}
 \le C\|h\|_{W^{2,4}}
 \le C\|H\|_{L^4}
 \le C\Lambda_1,
\end{equation*}
and $\nabla^2h=A-h\sigma$ then gives
$\|h\|_{C^2}\le C(\psi_-,\psi_+)$.

\emph{Step 2: positivity of the support function.}
By \eqref{eq:lowercurv}, $W[h-\kappa_1]\ge0$, so the
$1$-homogeneous extension of $h-\kappa_1$ is convex and
$h-\kappa_1$ is the support function of a nonempty compact
convex set $K$, possibly with empty interior \cite{schneider2014}.  The Steiner point map is
linear and vanishes on constants, since $\int_{\sph}\nu\dA=0$; by
\eqref{eq:perturbeddata} the Steiner point of $K$ is the origin.
The Steiner point of such a set lies in its relative interior
\cite{schneider2014}, so
$0\in K$ and the support function of $K$ is nonnegative:
\begin{equation*}
 h\ge\kappa_1.
\end{equation*}

\emph{Step 3: H\"older bounds.}
In a fixed finite atlas,
$W[h]_{ij}=\partial_i\partial_jh-\Gamma^k_{ij}\partial_kh+h\sigma_{ij}$,
so the equation reads $\sqrt{\det(W[h]_{ij})}=\sqrt{\psi\det\sigma}$.
On the corridor \eqref{eq:uppercurv}--\eqref{eq:lowercurv} it is
uniformly elliptic, and by Step~1 and \eqref{eq:databounds} its
lower-order terms and right-hand side are bounded in $C^{0,1}$ by
$C(\psi_-,\psi_+)$.  Such equations in two variables admit interior
$C^{2,\alpha}$ estimates without concavity
\cite[Section~17.3]{gilbargtrudinger2001}: there is
$\alpha_*=\alpha_*(\psi_-,\psi_+)\in(0,1)$ such that, for
$0<\alpha<\alpha_*$, finitely many charts give
\begin{equation}\label{eq:C2alpha}
 \|h\|_{C^{2,\alpha}}
 \le C(\psi_-,\psi_+,\alpha).
\end{equation}
Let $X_j(\nu)=\mathbf e_j\times\nu$, $j=1,2,3$, be the
rotation fields, so that $\nabla h=\sum_j(X_jh)\,X_j$.  Since $W$
commutes with rotations, differentiating $\det W[h]=\psi$ along their
flows gives
\begin{equation}\label{eq:killing}
 \mathcal L_{h}(X_j h)=X_j\psi,
 \qquad j=1,2,3.
\end{equation}
The coefficients of $\mathcal L_h$ are bounded in
$C^{0,\alpha}$ by \eqref{eq:C2alpha}, so the Schauder estimate
\cite[Section~6.1]{gilbargtrudinger2001} gives
$\|X_jh\|_{C^{2,\alpha}}\le C(\|X_j\psi\|_{C^{0,\alpha}}+\|X_jh\|_{C^0})
\le C(\|\psi\|_{C^{1,\alpha}}+\|h\|_{C^1})\le C$,
hence $\|h\|_{C^{3,\alpha}}\le C$.  The coefficients are now bounded in
$C^{1,\alpha}$; differentiating \eqref{eq:killing} once in each chart
and repeating the estimate, with the regularity theorem of
\cite[Section~6.4]{gilbargtrudinger2001}, gives
$\|X_jh\|_{C^{3,\alpha}}\le C$, hence
\begin{equation}\label{eq:C4alphabound}
 \|h\|_{C^{4,\alpha}}\le M_1,
\end{equation}
with $M_1$ depending only on $\psi_-,\psi_+,\alpha$.

\emph{Step 4: Sobolev bounds.}
By \eqref{eq:C4alphabound}, the linear estimate \eqref{eq:localtame}
applies to $X_jh$ in \eqref{eq:killing} and gains one derivative of $h$
at a time; iterating gives \eqref{eq:tame} (Appendix~\ref{app:four}).
\end{proof}

The space discretization enters the defect only through the
velocity solves: at a state $h\in V_L$, the computed velocity
$v=\mathcal V_L(h)(\rho-1)$ solves \eqref{eq:discretevelocity}, not the
linearized equation $\mathcal L_hv=\rho-1$.  Its failure to solve the
latter is the residual
\begin{align}\label{eq:zetadef}
\zeta[h]:=\mathcal{L}_{h}v-(\rho-1)=\mathcal L_h(v-\hat v),
\end{align}
the second form holding for $\hat v$ of
\eqref{eq:vhat}, which exists for $h$ near the exact path by
Proposition~\ref{prop:infsup}(ii).

\begin{proposition}[residual splitting and exact defect recursions]
\label{prop:discretedefect}
Let $h^n,v^n,k_1,k_2,g^n$ be as in Algorithm~\ref{alg:scheme}, assume
that every velocity solve below is well posed, and let $r_n$ be the
defect of \eqref{eq:exactperturbed}.  Then $r_0=0$, and:

(i) \emph{Residual splitting.}  For every state $h\in V_L$ at which
$\zeta[h]$ and $\hat v$ of \eqref{eq:vhat} are defined,
\begin{align}\label{eq:zetarep}
\zeta[h]=\bar a_L+(I-P_L)\,\mathcal{L}_{h}\bigl(v-\hat v\bigr),
\end{align}
where $\bar a_L$ is given by \eqref{eq:abar} and the second
term is orthogonal to $V_L$; in particular, $P_L\,\zeta[h]=\bar a_L$.

(ii) \emph{Euler variant.}  With $\zeta_n:=\zeta[h^n]$,
\begin{align}\label{eq:defecteulerL}
r_{n+1}=r_n+\tau\,\zeta_n+\tau^{2}\det W[v^n] .
\end{align}

(iii) \emph{Trapezoidal variant.}  With $\zeta_n^{\mathrm{st}}:=\zeta[g^n]$
the residual of the stage state,
\begin{align}\label{eq:defectheunL}
\begin{aligned}
r_{n+1}&=r_n+\tfrac{\tau}{2}\bigl(\zeta_n+\zeta_n^{\mathrm{st}}\bigr)
+\tfrac{\tau^{2}}{4}\det W[k_1-k_2],\\
b_{h^n}(k_1-k_2,q)&=\tau\,b_{k_1}(k_2,q)
\qquad\text{for all }q\in V_L .
\end{aligned}
\end{align}
\end{proposition}

\begin{proof}
\emph{Step 1: the residual splitting.}  Testing
\eqref{eq:discretevelocity} with $q\in V_L$ and using \eqref{eq:greens} and
\eqref{eq:loadrepresenter} give
\begin{align*}
(\zeta[h],q)=\sum_{k=1}^3\lambda_k(\nu_k,q)+(a_L,q)
\qquad\text{for all }q\in V_L ;
\end{align*}
by Corollary~\ref{cor:lambda},
$\sum_k\lambda_k(\nu_k,q)=(\lambda\cdot\nu,q)=-(P_{\mathcal N}a_L,q)$,
so the right-hand side equals
$(\bar a_L,q)$ by \eqref{eq:abar}, and hence $P_L\zeta[h]=\bar a_L$.  Now \eqref{eq:zetarep} follows from \eqref{eq:zetadef}.

\emph{Step 2: the Euler recursion.}  By
\eqref{eq:zetadef}, the computed velocity satisfies
$\mathcal{L}_{h^n}v^n=(\rho-1)+\zeta_n$.  The expansion
\eqref{eq:quadratic} applied to $h^{n+1}=h^n+\tau v^n$ and the path
identity $\rho_{t_{n+1}}=\rho_{t_n}+\tau(\rho-1)$ of
Lemma~\ref{lem:path} give
\begin{align*}
r_{n+1}
&=\det W[h^n]+\tau\,\mathcal{L}_{h^n}v^n+\tau^{2}\det W[v^n]
-\rho_{t_n}-\tau(\rho-1)\\
&=r_n+\tau\,\zeta_n+\tau^{2}\det W[v^n].
\end{align*}

\emph{Step 3: the trapezoidal recursion.}  The two stages solve
\begin{align*}
\mathcal{L}_{h^n}k_1&=(\rho-1)+\zeta_n,\\
\mathcal{L}_{h^n}k_2
&=\mathcal{L}_{g^n}k_2-\tau\,\mathcal{L}_{k_1}k_2
=(\rho-1)+\zeta_n^{\mathrm{st}}-2\tau D\bigl(W[k_1],W[k_2]\bigr),
\end{align*}
by \eqref{eq:zetadef},
$\mathcal{L}_{g^n}=\mathcal{L}_{h^n}+\tau\,\mathcal{L}_{k_1}$, and
\eqref{eq:Lop}.  The update
$h^{n+1}=h^n+\tfrac{\tau}{2}(k_1+k_2)$ in \eqref{eq:quadratic} and the
path identity give the first line below; adding the stage
equations gives
$\tfrac{\tau}{2}\mathcal{L}_{h^n}(k_1+k_2)=\tau(\rho-1)
+\tfrac{\tau}{2}(\zeta_n+\zeta_n^{\mathrm{st}})-\tau^2D(W[k_1],W[k_2])$,
hence the second:
\begin{align*}
r_{n+1}
&=r_n+\tfrac{\tau}{2}\,\mathcal{L}_{h^n}(k_1+k_2)
+\tfrac{\tau^{2}}{4}\det W[k_1+k_2]-\tau(\rho-1)\\
&=r_n+\tfrac{\tau}{2}\bigl(\zeta_n+\zeta_n^{\mathrm{st}}\bigr)
+\tfrac{\tau^{2}}{4}\Bigl(\det W[k_1+k_2]
-4D\bigl(W[k_1],W[k_2]\bigr)\Bigr).
\end{align*}
By \eqref{eq:cofactorids} and \eqref{eq:mixeddisc},
$\det(A+B)-4D(A,B)=\det(A-B)$ for symmetric $A,B$; with $A=W[k_1]$ and
$B=W[k_2]$, the last bracket is $\det W[k_1-k_2]$.  For the second line
of \eqref{eq:defectheunL}, both stage systems have the same data term
and, by Corollary~\ref{cor:lambda}, the same multiplier, so subtracting
them gives $b_{h^n}(k_1,q)=b_{g^n}(k_2,q)$ for $q\in V_L$.  Since
$g^n=h^n+\tau k_1$ and $h\mapsto b_h$ is linear, the right-hand side is
$b_{h^n}(k_2,q)+\tau\,b_{k_1}(k_2,q)$, which gives the second line of
\eqref{eq:defectheunL}.
\end{proof}

Set $\Theta(h):=1+\|h\|_{H^{s+2}}$ and
$\varepsilon_2:=\tfrac18\min\{\varepsilon_0,\varepsilon_1,\beta/4,\kappa_0\}$,
with $\varepsilon_0$ and $L_0$ from Lemma~\ref{lem:discreteinfsup},
$\varepsilon_1$ from Lemma~\ref{lem:stationary}, and $M_1$ from
Lemma~\ref{lem:apriori} with $(\psi_-,\psi_+)=(\delta,\Lambda_*)$.  The
iterates satisfy \eqref{eq:zetahyp} with $\varepsilon_2$ and $M_1$ in
place of $2\varepsilon_2$ and $M_1+1$; the extra room is for the
trapezoidal stage state $g^n=h^n+\tau k_1$.

\begin{lemma}[residual and velocity bounds]\label{lem:zeta}
Let $L\ge L_0$ and let $h\in V_L^\perp$ satisfy
\begin{align}\label{eq:zetahyp}
\|W[h]-W[h(t)]\|_{L^\infty}\le2\varepsilon_2
\;\text{ for some }t\in[0,1],
\qquad
\|h\|_{C^{4,\alpha}}\le M_1+1 .
\end{align}
Then $v=\mathcal{V}_L(h)(\rho-1)$ and $\hat v$ of \eqref{eq:vhat}
are well defined.  With $\Theta=\Theta(h)$ and generic constants $C$ as
in Section~\ref{sec:analysis}, which depend on $\varepsilon_2$ and $M_1$
but not on $h$, $L$, or $\tau$, the function $\hat v$ and the Galerkin
gap $v-P_L\hat v$ satisfy
\begin{align}\label{eq:whatbounds}
\|\hat v\|_{C^{2,\alpha}}\le C,
\qquad
\|\hat v\|_{H^{s+2}}\le C\,\Theta,
\qquad
\|v-P_L\hat v\|_{H^1}\le C\,\Theta\,L^{-(s+1)} ,
\end{align}
the computed velocity $v$ satisfies
\begin{align}\label{eq:vbounds}
\|v\|_{H^1}\le C,
\qquad
\|v\|_{H^{s+2}}\le C\,\Theta,
\qquad
\|W[v]\|_{L^\infty}\le C+C\,\Theta\,L^{1+\alpha-s} ,
\end{align}
and the residual $\zeta[h]$ of \eqref{eq:zetadef} satisfies
\begin{align}
&\|\zeta[h]\|_{(\Hperp)'}\le C\,\Theta\,L^{-(s+1)},
\qquad
\|\zeta[h]\|_{C^{2,\alpha}}\le C\,\Theta^{2}L^{3+2\alpha-s},
\label{eq:zetadual}\\
&\|\zeta[h]\|_{H^{s}}\le C\,\Theta+C\,\Theta^{2}L^{1+2\alpha-s} .
\label{eq:zetaHs}
\end{align}
\end{lemma}

\begin{proof}
Throughout, $\eta:=v-P_L\hat v\in V_L^\perp$; every inverse
estimate acts on $\eta$ alone.

\emph{Step 1: well-posedness.}  Since
$2\varepsilon_2\le\min\{\varepsilon_0,\kappa_0/4,\beta/16\}$, the tube
bound in \eqref{eq:zetahyp} gives the hypothesis of
Lemma~\ref{lem:discreteinfsup}, places $h$ in
$\mathcal{W}(\kappa_0/2,2\Lambda_0)$, and, by
Proposition~\ref{prop:infsup}(ii), gives the inf-sup bound at $h$ with
$\beta/2$ in place of $\beta$.  Hence $v$ and $\hat v$ exist with
$\|v\|_{H^1}+\|\hat v\|_{H^1}\le C$, and Lemma~\ref{lem:lineartame}
with $f=\rho-1$ and, by the second bound in \eqref{eq:zetahyp},
$M_0=M_1+1$ gives the bounds on $\hat v$ in \eqref{eq:whatbounds}.

\emph{Step 2: the Galerkin gap.}  By \eqref{eq:greens},
$b_{h}(\hat v,q)=-(\rho-1,q)$ for $q\in H^1$.  Testing
\eqref{eq:discretevelocity} with $q\in V_L^\perp$ removes the
multiplier, and \eqref{eq:loadrepresenter} with \eqref{eq:abar} gives
$b_h(v,q)=-(\rho-1,q)-(\bar a_L,q)$.  Subtracting, and then using
\eqref{eq:discreteinfsup}, \eqref{eq:approx} with Step~1, and
\eqref{eq:loadbounds},
\begin{align*}
b_h(\eta,q)&=b_h(\hat v-P_L\hat v,q)-(\bar a_L,q),
\qquad q\in V_L^\perp,\\
\|\eta\|_{H^1}
&\le C\bigl(\|\hat v-P_L\hat v\|_{H^1}
+\|\bar a_L\|_{(\Hperp)'}\bigr)
\le C\,\Theta\,L^{-(s+1)} .
\end{align*}

\emph{Step 3: velocity bounds.}  Since
$W[v]=W[\eta]+W[\hat v]-W[(I-P_L)\hat v]$, Steps~1--2,
\eqref{eq:approx}, \eqref{eq:inverse},
Lemma~\ref{lem:Winfty}, and
$H^{3+\alpha}(\sph)\hookrightarrow W^{2,\infty}(\sph)$ give
\begin{align*}
\|v\|_{H^{s+2}}
&\le\|P_L\hat v\|_{H^{s+2}}+CL^{s+1}\|\eta\|_{H^1}
\le C\,\Theta,\\
\|W[v]\|_{L^\infty}
&\le CL^{2}\|\eta\|_{H^1}
+\|W[(I-P_L)\hat v]\|_{L^\infty}+\|W[\hat v]\|_{L^\infty}
\le C\,\Theta\,L^{1+\alpha-s}+C .
\end{align*}

\emph{Step 4: residual bounds.}  By
Proposition~\ref{prop:discretedefect}(i), the residual splits as
$\zeta[h]=\bar a_L+\zeta^\perp$, where
$\zeta^\perp:=(I-P_L)\mathcal L_h(v-\hat v)$ is orthogonal to $V_L$; by
\eqref{eq:loadbounds}, $\bar a_L$ satisfies
\eqref{eq:zetadual}--\eqref{eq:zetaHs}, so it suffices to bound
$\zeta^\perp$.  Since
$v-\hat v=\eta-(I-P_L)\hat v$, Steps~1--2,
\eqref{eq:approx}, \eqref{eq:inverse}, and
$H^{3+2\alpha}\hookrightarrow C^{2,\alpha}$ give
\begin{align*}
\|v-\hat v\|_{H^2}\le C\,\Theta\,L^{-s},
\qquad
\|v-\hat v\|_{H^{5+2\alpha}}\le C\,\Theta\,L^{3+2\alpha-s},
\qquad
\|v-\hat v\|_{C^{2,\alpha}}\le C\,\Theta\,L^{1+2\alpha-s}.
\end{align*}
By \eqref{eq:Lop},
$\mathcal L_h(v-\hat v)=C[h]^{ij}W[v-\hat v]_{ij}$, where
$\|C[h]\|_{C^{\alpha}}\le C\|h\|_{C^{4,\alpha}}\le C$ by
\eqref{eq:zetahyp} and, since $5+2\alpha\le s+2$,
$\|C[h]\|_{H^{3+2\alpha}}+\|C[h]\|_{H^{s}}\le C\,\Theta$.  Hence
\eqref{eq:moser} and Steps~1 and~3 give
\begin{align*}
\|\zeta^\perp\|_{C^{2,\alpha}}
&\le C\|\zeta^\perp\|_{H^{3+2\alpha}}
\le C\,\Theta^{2}L^{3+2\alpha-s},\\
\|\zeta^\perp\|_{H^{s}}
&\le C\bigl(\Theta\,\|v-\hat v\|_{C^{2,\alpha}}
+\|v-\hat v\|_{H^{s+2}}\bigr)
\le C\,\Theta+C\,\Theta^{2}L^{1+2\alpha-s} .
\end{align*}
In the dual norm, $\zeta^\perp\perp V_L$ gives
$(\zeta^\perp,q)=(\zeta^\perp,(I-P_L)q)$ for $q\in\Hperp$, and
\eqref{eq:approx} bounds $\|(I-P_L)q\|_{L^2}$ by $CL^{-1}\|q\|_{H^1}$.
With $\|\zeta^\perp\|_{L^2}\le C\|v-\hat v\|_{H^2}$, this extra factor
$L^{-1}$ gives $\|\zeta^\perp\|_{(\Hperp)'}\le C\,\Theta\,L^{-(s+1)}$.
\end{proof}

\subsection{Convergence of the fully discrete continuation}

Both proofs use induction to preserve convexity, the defect
bounds, and $C^2$-closeness to the exact path.
Starting from \eqref{eq:deltalambda}--\eqref{eq:alphachoice},
we fix the constants in the order
\[
\begin{aligned}
 (\delta,\gamma,\Lambda_*)&\longrightarrow\alpha
 \longrightarrow(\kappa_1,\Lambda_1,M_1)
 \longrightarrow\varepsilon_2\\
 &\longrightarrow C_\star
 \longrightarrow B
 \longrightarrow(L_0,\tau_0),
\end{aligned}
\]
each from the data and the constants before it; $C_\star$ is defined in
\eqref{eq:Cstar} below, and $L_0,\tau_0$ in the two proofs.  Neither of
those two thresholds involves the other parameter, so no condition
couples $\tau$ and $L$.

As noted before Lemma~\ref{lem:apriori}, for a normalized,
strictly convex iterate $h^n$ with $\|r_n\|_{C^{2,\alpha}}\le\delta$,
Lemma~\ref{lem:apriori} gives the uniform bounds
\eqref{eq:corridorA} and, by \eqref{eq:tame} at order $s$,
\[
 \Theta(h^n)\le C\bigl(1+\|r_n\|_{H^s}\bigr).
\]
For the trapezoidal stage $g^n=h^n+\tau k_1$,
\eqref{eq:vbounds} also gives
\begin{align}\label{eq:stagetheta}
 \Theta(g^n)\le(1+C\tau)\Theta(h^n)\le2\Theta(h^n)
\end{align}
for sufficiently small $\tau$, independently of $L$.
Thus, once the hypotheses of Lemma~\ref{lem:zeta} are
verified at each stage, \eqref{eq:zetaHs} yields
\begin{equation}\label{eq:Cstar}
 \|\zeta[h]\|_{H^s}
 \le C_\star\bigl(1+\|r_n\|_{H^s}\bigr)
      +C\Theta(h)^2L^{1+2\alpha-s},
\end{equation}
for $h=h^n$ and, in the trapezoidal scheme, $h=g^n$.
Like the constants $C$ above, $C_\star$ depends only on the
data and on what is already fixed, never on $B$, $L$, $\tau$, or $n$.
We therefore choose
\(
 B:=2e^{C_\star}(2+C_\star),
\)
so that the Gronwall estimate below gives the stronger bound
$\|r_n\|_{H^s}\le B/2$.  In the two proofs, arguments in
parentheses list the constants a quantity depends on, as in $C(B)$ and
$\tau_0(B,\delta)$.

\begin{theorem}[convergence]\label{thm:main}
Let $\rho$ satisfy the hypotheses stated at the beginning of this section.
There exist $\tau_0>0$, $L_0\in\mathbb{N}$, and $C>0$ such that,
for every $N\in\mathbb N$ with $\tau=1/N\le\tau_0$ and every integer
$L\ge L_0$, the Euler
variant of Algorithm~\ref{alg:scheme} with $M=0$ is well defined, all its
states satisfy
$W[h^n]\ge\tfrac{\kappa_0}{2}\sigma$, and
\begin{align}\label{eq:mainestimate}
\max_{0\le n\le N}\,\bigl\|h^n-h(t_n)\bigr\|_{H^1(\sph)}
\;\le\;C\bigl(\tau+L^{-(s+1)}\bigr).
\end{align}
\end{theorem}

\begin{proof}
We show by induction that $W[h^m]>0$ and
\begin{align}\label{eq:finduction}
\begin{split}
(\mathrm D_\alpha)\quad&\|r_m\|_{C^{2,\alpha}}\le\delta,\\
(\mathrm D_s)\quad&\|r_m\|_{H^{s}}\le B,\\
(\mathrm T)\quad&\|h^m-h(t_m)\|_{C^{2}}\le\tfrac12\varepsilon_2
\end{split}
\end{align}
hold for $0\le m\le n$.  $(\mathrm D_\alpha)$ keeps the
data admissible for Lemma~\ref{lem:apriori}, $(\mathrm D_s)$ bounds
$\|h^m\|_{H^{s+2}}$, and $(\mathrm T)$ keeps $h^m$ in the tube of
Lemma~\ref{lem:stationary}.  All three hold at $m=0$, where
$h^0=1=h(0)$ and $r_0=0$.

\emph{Step 1: what the hypotheses give.}  Each $h^m$ is strictly
convex and normalized, so, as noted at the beginning of this
subsection, $(\mathrm D_\alpha)$ and Lemma~\ref{lem:apriori} give
\eqref{eq:corridorA} and, with
$(\mathrm D_s)$,
\begin{align}\label{eq:hsbound}
\Theta(h^m)\le C\bigl(1+\|r_m\|_{H^s}\bigr)\le C(B).
\end{align}

\emph{Step 2: Lemma~\ref{lem:zeta} applies at each $h^m$.}  Its first
hypothesis follows from $(\mathrm T)$ and the linearity of $W$,
\begin{align*}
\|W[h^{m}]-W[h(t_{m})]\|_{L^\infty}\le2\|h^{m}-h(t_{m})\|_{C^{2}}\le\varepsilon_2,
\end{align*}
and its second from \eqref{eq:corridorA}.  So
\eqref{eq:whatbounds}--\eqref{eq:zetaHs} hold at $h^m$ for
$m\le n$, with $\Theta=\Theta(h^{m})$.

\emph{Step 3: $(\mathrm D_\alpha)$ and $(\mathrm D_s)$ at level $n+1$.}  By
\eqref{eq:vbounds}, \eqref{eq:hsbound}, and
$H^{s+2}(\sph)\hookrightarrow C^{4,\alpha}(\sph)$,
\begin{align*}
\|v^n\|_{H^{s+2}}\le C(B),
\qquad
\|W[v^n]\|_{H^{s}}+\|W[v^n]\|_{C^{2,\alpha}}\le C\|v^n\|_{H^{s+2}} .
\end{align*}
Since $\det W[v^n]$ is quadratic in the entries of $W[v^n]$,
\eqref{eq:moser} and the product rule give
\begin{align*}
\|\det W[v^n]\|_{C^{2,\alpha}}+\|\det W[v^n]\|_{H^{s}}\le C(B).
\end{align*}
Summing \eqref{eq:defecteulerL} from $r_0=0$,
\begin{align*}
\|r_{n+1}\|_{C^{2,\alpha}}
\le\sum_{m\le n}\tau\bigl(\|\zeta_m\|_{C^{2,\alpha}}
+\tau\,\|\det W[v^m]\|_{C^{2,\alpha}}\bigr)
\le C(B)\bigl(L^{3+2\alpha-s}+\tau\bigr)
\le\delta
\end{align*}
for $L\ge L_0(B,\delta)$ and $\tau\le\tau_0(B,\delta)$, since
$3+2\alpha-s<0$ by \eqref{eq:alphachoice} and $s\ge4$;
this is $(\mathrm D_\alpha)$.  By \eqref{eq:Cstar} and \eqref{eq:hsbound} at $h^n$, once
$C(B)L^{1+2\alpha-s}\le\tfrac12$ and $\tau\,C(B)\le\tfrac12$, the
$\Theta^{2}$ term in \eqref{eq:Cstar} and $\tau^{2}\det W[v^n]$ each add
at most $\tfrac{\tau}{2}$ to \eqref{eq:defecteulerL} per step, so
\begin{align*}
\|r_{n+1}\|_{H^{s}}
\le\|r_n\|_{H^{s}}
+\tau\,C_\star\bigl(1+\|r_n\|_{H^{s}}\bigr)+\tau ,
\end{align*}
and the discrete Gronwall inequality, with $r_0=0$ and $N\tau=1$,
gives $\|r_{n+1}\|_{H^{s}}\le e^{C_\star}(2+C_\star)=\tfrac{B}{2}$;
this is $(\mathrm D_s)$.

\emph{Step 4: convexity and the corridor at level $n+1$.}  By
\eqref{eq:vbounds} and \eqref{eq:hsbound},
\begin{align*}
\|W[h^{n+1}]-W[h^n]\|_{L^\infty}
=\tau\,\|W[v^n]\|_{L^\infty}
\le C(B)\,\tau
\le\tfrac{\kappa_1}{2}
\end{align*}
if $\tau$ is small enough; so
$W[h^{n+1}]\ge\tfrac{\kappa_1}{2}\sigma>0$, and Step~1 now applies at
level $n+1$ and restores \eqref{eq:corridorA} with the same constants.

\emph{Step 5: $(\mathrm T)$ at level $n+1$ and the error bound.}  By
\eqref{eq:Lhh} and \eqref{eq:greens}, $2(\det W[v^m],q)=-b_{v^m}(v^m,q)$
for $q\in\Hperp$; since $C[v^m]$ and $\tr W[v^m]$ are bounded by
$2\|W[v^m]\|_{L^\infty}$, Cauchy--Schwarz and \eqref{eq:vbounds} give
\begin{align*}
\|\det W[v^m]\|_{(\Hperp)'}
\le C\,\|W[v^m]\|_{L^\infty}\|v^m\|_{H^1}\le C(B).
\end{align*}
Summing \eqref{eq:defecteulerL} in the dual norm, with
\eqref{eq:zetadual} and $N\tau=1$,
\begin{align}\label{eq:rdual}
\|r_{n+1}\|_{(\Hperp)'}
\le\sum_{m\le n}\tau\bigl(\|\zeta_m\|_{(\Hperp)'}
+\tau\,\|\det W[v^m]\|_{(\Hperp)'}\bigr)
\le C(B)\bigl(\tau+L^{-(s+1)}\bigr).
\end{align}
By Step~3, \eqref{eq:pathregularity}, and
$H^{s+2}(\sph)\hookrightarrow C^{2}(\sph)$, one step moves $h^n$ and
$h(t_n)$ by at most $C(B)\tau$ in $C^2$, so $(\mathrm T)$ gives
\begin{align*}
\|h^{n+1}-h(t_{n+1})\|_{C^2}
\le\|h^n-h(t_n)\|_{C^2}+C(B)\,\tau
\le\varepsilon_2
\end{align*}
once $C(B)\,\tau\le\tfrac12\varepsilon_2$.  Hence
$\|W[h^{n+1}]-W[h(t_{n+1})]\|_{L^\infty}\le2\varepsilon_2\le\varepsilon_1$,
and Lemma~\ref{lem:stationary} with \eqref{eq:rdual} gives
\begin{align}\label{eq:h1error}
\|h^{n+1}-h(t_{n+1})\|_{H^1}
\le C\,\|r_{n+1}\|_{(\Hperp)'}
\le C(B)\bigl(\tau+L^{-(s+1)}\bigr).
\end{align}
By Steps~3 and~4 and \eqref{eq:pathregularity},
$\|h^{n+1}-h(t_{n+1})\|_{H^{s+2}}\le C(B)$.  Interpolating this with
\eqref{eq:h1error} in \eqref{eq:sobolev} at
$\theta=\tfrac{s-1-\alpha}{s+1}$, and using
$H^{3+\alpha}\hookrightarrow C^{2}$,
\begin{align*}
\|h^{n+1}-h(t_{n+1})\|_{C^{2}}
\le C(B)\bigl(\tau+L^{-(s+1)}\bigr)^{\theta}
\le\tfrac12\varepsilon_2
\end{align*}
once $\tau\le\tau_0$ and $L\ge L_0$, fixed from $B$ and $\varepsilon_2$
and not from each other.  This is $(\mathrm T)$ at level $n+1$, so the
induction closes, and the maximum of \eqref{eq:h1error} over $n$ is
\eqref{eq:mainestimate}.  Finally,
$(\mathrm T)$ keeps $W[h^n]$ within $\varepsilon_2$ of
$W[h(t_n)]$, so \eqref{eq:uniformbounds} and
$\varepsilon_2\le\tfrac{\kappa_0}{8}$ give
$W[h^n]\ge(\kappa_0-\varepsilon_2)\,\sigma\ge\tfrac{\kappa_0}{2}\,\sigma$.
\end{proof}

For a generic ordinary differential equation the explicit
trapezoidal method is second order.  Here the determinant is exactly quadratic in $h$ because $W[h]$
is a $2\times2$ matrix, and the data increment $\tau(\rho-1)$ is exact because
$\rho_t$ is affine in $t$.  Proposition~\ref{prop:discretedefect}
therefore expands one step with no remainder, leaving the single
quadratic term $\tfrac{\tau^2}{4}\det W[k_1-k_2]$; since
$k_1-k_2=O(\tau)$ by the second line of \eqref{eq:defectheunL}, it is
$O(\tau^4)$, and the accumulated defect is third order, as measured in
Section~\ref{subsec:expA}.  If either structure is absent, the trapezoidal
rule generally falls back to its classical second order.

\begin{theorem}[third-order convergence of the explicit trapezoidal variant]
\label{thm:trap}
Let $\rho$ satisfy the hypotheses stated at the beginning of this section.
There exist $\tau_0>0$, $L_0\in\mathbb N$, and $C>0$ such that,
for every $N\in\mathbb N$ with $\tau=1/N\le\tau_0$ and every integer
$L\ge L_0$, the trapezoidal variant of Algorithm~\ref{alg:scheme}
with $M=0$ is well defined, all its iterates and stage states satisfy
\(
W[h^n],W[g^n]\ge\tfrac{\kappa_0}{2}\sigma,
\)
and
\begin{align}\label{eq:trapestimate}
\max_{0\le n\le N}
\|h^n-h(t_n)\|_{H^1(\sph)}
\le C\bigl(\tau^3+L^{-(s+1)}\bigr).
\end{align}
\end{theorem}

\begin{proof}
We repeat the induction \eqref{eq:finduction}, replacing $(\mathrm T)$ by
\[
(\mathrm T_{\mathrm{tr}})\qquad
\|h^m-h(t_m)\|_{C^2}\le\tfrac14\varepsilon_2 .
\]
The smaller tube leaves room for the stage state; assume the
induction through level $n$.

\emph{Step 1: what the hypotheses give.}
As in Theorem~\ref{thm:main},
\eqref{eq:corridorA} and \eqref{eq:hsbound} hold at $h^n$.

\emph{Step 2: Lemma~\ref{lem:zeta} applies at $h^n$ and at $g^n=h^n+\tau k_1$.}
At $h^n$ it applies as in Step~2 of the proof of
Theorem~\ref{thm:main}, since
$(\mathrm T_{\mathrm{tr}})$ implies $(\mathrm T)$; by \eqref{eq:vbounds},
$\|k_1\|_{C^{4,\alpha}}\le C\|k_1\|_{H^{s+2}}\le C(B)$.

At $g^n$, this bound, $(\mathrm T_{\mathrm{tr}})$, and \eqref{eq:corridorA}
give, after reducing $\tau_0(B)$ independently of $L$,
\begin{align}\label{eq:stagebounds}
\|W[g^n]-W[h(t_n)]\|_{L^\infty}
\le{2\cdot\tfrac14}\varepsilon_2+C(B)\tau\le\varepsilon_2,
\quad
\|g^n\|_{C^{4,\alpha}}\le M_1+C(B)\tau\le M_1+1,
\end{align}
and \eqref{eq:stagetheta} with \eqref{eq:hsbound} gives
$\Theta(g^n)\le2\Theta(h^n)\le C(B)$.
Hence Lemma~\ref{lem:zeta} holds at both stages with $\Theta\le C(B)$;
in particular $\|k_j\|_{H^1}\le C$ and $\|k_j\|_{H^{s+2}}\le C(B)$ for
$j=1,2$.

\emph{Step 3: $(\mathrm D_\alpha)$ and $(\mathrm D_s)$ at level $n+1$.}
We follow Step~3 of the proof of Theorem~\ref{thm:main},
with $k_1-k_2$ in place of $v^n$.  By Step~2,
$\|k_1-k_2\|_{H^{s+2}}\le C(B)$, so $\det W[k_1-k_2]$ is
bounded by $C(B)$ in $C^{2,\alpha}$ and in $H^s$.  Summing the first line of \eqref{eq:defectheunL}, with
\eqref{eq:zetadual} at every stage,
\[
\|r_{n+1}\|_{C^{2,\alpha}}
\le\sum\nolimits_{m\le n}\bigl(\tau\,C(B)L^{3+2\alpha-s}+\tau^2C(B)\bigr)
\le C(B)\bigl(L^{3+2\alpha-s}+\tau\bigr)\le\delta
\]
for $L\ge L_0(B,\delta)$ and $\tau\le\tau_0(B,\delta)$, chosen
independently; this is $(\mathrm D_\alpha)$.  By Step~2,
\eqref{eq:Cstar} holds at both stages, so the $H^s$ recursion of
Theorem~\ref{thm:main} carries over, after enlarging $L_0$ and reducing $\tau_0$
separately, and gives $(\mathrm D_s)$:
$\|r_{n+1}\|_{H^s}\le e^{C_\star}(2+C_\star)=\tfrac12B$.

\emph{Step 4: convexity and the corridor at level $n+1$.}
Since $h^{n+1}-h^n=\tfrac\tau2(k_1+k_2)$, Step~2 gives
$\|W[h^{n+1}]-W[h^n]\|_{L^\infty}\le C(B)\tau$, as in Step~4 of the
proof of Theorem~\ref{thm:main}.  So $h^{n+1}$ is strictly convex, and with $(\mathrm D_\alpha)$
from Step~3, Lemma~\ref{lem:apriori} restores \eqref{eq:corridorA} at
level $n+1$.

\emph{Step 5: the stage difference.}
The second line of \eqref{eq:defectheunL} reads
$b_{h^n}(k_1-k_2,q)=\tau\,b_{k_1}(k_2,q)$, $q\in V_L$; the inf-sup
bound \eqref{eq:discreteinfsup}, \eqref{eq:bsym}--\eqref{eq:bbound},
and Step~2 give
\begin{align}\label{eq:k1k2H1}
\|k_1-k_2\|_{H^1}
\le C\tau\|k_1\|_{C^2}\|k_2\|_{H^1}
\le C(B)\tau .
\end{align}

For $W[k_1-k_2]$ in $L^\infty$, let
$\hat v_1,\hat v_2\perp\mathcal N$ solve \eqref{eq:vhat} at $h^n$
and at $g^n$; since
$\mathcal L_{g^n}=\mathcal L_{h^n}+\tau\mathcal L_{k_1}$, the
difference of the $\hat v_j$ solves
\[
\mathcal L_{h^n}\hat v_{12}=2\tau D\bigl(W[k_1],W[\hat v_2]\bigr),
\qquad
\hat v_{12}:=\hat v_1-\hat v_2\perp\mathcal N .
\]
By \eqref{eq:greens} and Step~2, the source is orthogonal to
$\mathcal N$ and is $C(B)\tau$ in $C^\alpha$ and in $H^s$, so
Lemma~\ref{lem:lineartame} gives
$\|\hat v_{12}\|_{C^{2,\alpha}}+\|\hat v_{12}\|_{H^{s+2}}\le C(B)\tau$.
Now write
$k_1-k_2=(k_1-P_L\hat v_1)-(k_2-P_L\hat v_2)-(I-P_L)\hat v_{12}
+\hat v_{12}$, and bound $W[\cdot]$ of the first two pieces by
Lemma~\ref{lem:Winfty} with \eqref{eq:whatbounds}, of the third by
\eqref{eq:approx} with $H^{3+\alpha}\hookrightarrow W^{2,\infty}$, and
of the last by the bound on $\hat v_{12}$:
\begin{align}\label{eq:Wk1k2}
\|W[k_1-k_2]\|_{L^\infty}
\le C(B)\bigl(L^{2}L^{-(s+1)}+\tau L^{1+\alpha-s}+\tau\bigr)
\le C(B)\bigl(L^{1-s}+\tau\bigr).
\end{align}

\emph{Step 6: third order and the tube at level $n+1$.}
The dual bound of Step~5 of the proof of Theorem~\ref{thm:main}
holds for any element of $V_L$; with \eqref{eq:k1k2H1} and \eqref{eq:Wk1k2} it gives
\[
\|\det W[k_1-k_2]\|_{(\Hperp)'}
\le{C\|W[k_1-k_2]\|_{L^\infty}\|k_1-k_2\|_{H^1}}
\le C(B)\bigl(\tau^2+\tau L^{1-s}\bigr).
\]
Summing the first line of \eqref{eq:defectheunL} in the dual
norm, with \eqref{eq:zetadual} at every stage and $N\tau=1$,
\[
\|r_{n+1}\|_{(\Hperp)'}
\le C(B)\bigl(L^{-(s+1)}+\tau^3+\tau^2L^{1-s}\bigr)
\le C(B)\bigl(\tau^3+L^{-(s+1)}\bigr),
\]
the last step by Young's inequality
$\tau^2L^{1-s}\le\tfrac23\tau^3+\tfrac13L^{3(1-s)}$
and $3(1-s)\le-(s+1)$ for $s\ge2$.

As in Step~5 of the proof of Theorem~\ref{thm:main}, $(\mathrm T_{\mathrm{tr}})$, Step~2, and
\eqref{eq:pathregularity} give
$\|W[h^{n+1}]-W[h(t_{n+1})]\|_{L^\infty}
\le2\bigl(\tfrac14\varepsilon_2+C(B)\tau\bigr)\le\varepsilon_2\le\varepsilon_1$,
and Lemma~\ref{lem:stationary} yields
\begin{align}\label{eq:traph1}
\|h^{n+1}-h(t_{n+1})\|_{H^1}
\le C(B)\bigl(\tau^3+L^{-(s+1)}\bigr).
\end{align}
Interpolating \eqref{eq:traph1} against the $H^{s+2}$ bound
from Steps~3--4 and \eqref{eq:pathregularity} restores
$(\mathrm T_{\mathrm{tr}})$ after tightening $L_0$ and $\tau_0$
independently; the induction closes, and the maximum of
\eqref{eq:traph1} over $n$ is \eqref{eq:trapestimate}.  Finally
$(\mathrm T_{\mathrm{tr}})$ and \eqref{eq:stagebounds} keep $W[h^n]$ and
$W[g^n]$ within $\varepsilon_2$ of $W[h(t_n)]$, so
\eqref{eq:uniformbounds} and $\varepsilon_2\le\kappa_0/8$ give the
convexity stated for both families.
\end{proof}
\subsection{The Newton phase}\label{subsec:newton}

Theorems~\ref{thm:main} and~\ref{thm:trap} treat $M=0$.  The
Newton steps ($M>0$) converge instead to the solution of the discrete
problem, in which $a_L$ of \eqref{eq:loadrepresenter} represents the
quadrature error:
\begin{align}\label{eq:galerkinproblem}
 h_L^\ast\in V_L^\perp,\qquad
 (\det W[h_L^\ast]-\rho,q)=(a_L,q)
 \quad\text{for all }q\in V_L^\perp .
\end{align}

\begin{proposition}[the Newton phase]\label{prop:newton}
Assume the hypotheses of Theorem~\ref{thm:main}.  There exist
$\delta_1>0$ and $C>0$, and $L_1\ge L_0$ depending on $\delta_1$,
such that, for every $L\ge L_1$,
\eqref{eq:galerkinproblem} has exactly one solution in
\begin{align*}
 \mathbb B_L:=
 \left\{h\in V_L^\perp:
 \|h-P_Lh(1)\|_{H^1}\le\delta_1L^{-2}\right\}.
\end{align*}
This solution satisfies
\begin{align}\label{eq:galerkinerror}
 \|h_L^\ast-h(1)\|_{H^1}\le C L^{-(s+1)},
 \qquad
 W[h_L^\ast]\ge\tfrac{\kappa_0}{2}\sigma .
\end{align}
For every $h\in\mathbb B_L$, Lemma~\ref{lem:discreteinfsup} applies,
so the Newton step
$h^+=h+\mathcal V_L(h)\bigl(\rho-\det W[h]\bigr)$
is well defined and satisfies
\begin{align}\label{eq:newtonquad}
 \|h^+-h_L^\ast\|_{H^1}
 \le C L^2\|h-h_L^\ast\|_{H^1}^2 .
\end{align}
For sufficiently large $L$ and a suitable $c>0$, the continuation
endpoint $h^N$ lies in $\mathbb B_L$ provided
\[
 \tau\le cL^{-2}\quad\text{for Euler},
 \qquad
 \tau\le cL^{-2/3}\quad\text{for trapezoidal}.
\]
The Newton iterates then converge quadratically to $h_L^\ast$,
whose error is bounded in \eqref{eq:galerkinerror}.
\end{proposition}

\begin{proof}
Write $x_0:=P_Lh(1)$, $\eta:=(I-P_L)h(1)$, and $R:=\delta_1L^{-2}$, so
that $\mathbb B_L$ is the $H^1$-ball of radius $R$ about $x_0$.

 We first fix $\delta_1$
small (Steps~2 and~4), then $L_1$ large (Steps~2--4).

\emph{Step 1: the Newton map.}
Define $G_L:V_L^\perp\to(V_L^\perp)'$ by
$G_L(h)q=(\det W[h]-\rho-a_L,q)$, so that
\eqref{eq:galerkinproblem} reads $G_L(h_L^\ast)=0$.  By
\eqref{eq:quadratic} and \eqref{eq:greens}, with no remainder,
\[
 G_L(h+w)=G_L(h)+G_L'(h)w+(\det W[w],\cdot),
 \qquad
 G_L'(h)w=-b_h(w,\cdot).
\]
Testing the Newton step with $q\in V_L^\perp$ removes the multiplier
term.  Since $Q_L$ is exact on $q$ and on $\det W[h]\,q$
(Lemma~\ref{lem:Winfty} and \eqref{eq:QLexact}), \eqref{eq:loadrepresenter} gives
$Q_L\bigl((\rho-\det W[h])q\bigr)=-G_L(h)q$.  Hence the step is
Newton's method for $G_L$: $G_L'(h)(h^+-h)=-G_L(h)$.

\emph{Step 2: the bounds used in Steps~3 and~4.}
Since $G_L'(h)w=-b_h(w,\cdot)$ is linear in $h$ and
$(\det W[w],\cdot)=-\frac12b_w(w,\cdot)$ by \eqref{eq:greens} and
\eqref{eq:Lhh}, one bound on $b_z(w,q)$ controls both the change of
$G_L'$ and the quadratic term of Step~1: as in Step~5 of the proof of
Theorem~\ref{thm:main},
$|b_z(w,q)|\le2\|W[z]\|_{L^\infty}\|w\|_{H^1}\|q\|_{H^1}$, and by
\eqref{eq:Winfty} this is ${}\le CL^2\|z\|_{H^1}\|w\|_{H^1}\|q\|_{H^1}$ if $z\in V_L$.
Hence, for $h,z,w\in V_L^\perp$,
\begin{align}\label{eq:newtonlip}
 \|G_L'(h+z)-G_L'(h)\|_{V_L^\perp\to(V_L^\perp)'}\le CL^2\|z\|_{H^1},
 \qquad
 \|\det W[w]\|_{(V_L^\perp)'}\le CL^2\|w\|_{H^1}^2 .
\end{align}
For $h\in\mathbb B_L$, $h-h(1)=(h-x_0)-\eta$, so \eqref{eq:Winfty}, \eqref{eq:approx}, and \eqref{eq:pathregularity} give
\[
 \|W[h]-W[h(1)]\|_{L^\infty}
 \le CL^2\|h-x_0\|_{H^1}+\|W[\eta]\|_{L^\infty}
 \le C\delta_1+CL^{1+\alpha-s}
 \le\min\{\varepsilon_0,\kappa_0/2\},
\]
for $\delta_1$ small and $L\ge L_1$; so Lemma~\ref{lem:discreteinfsup} and \eqref{eq:uniformbounds} give,
uniformly for $h\in\mathbb B_L$,
\begin{align}\label{eq:newtoninverse}
 \|G_L'(h)^{-1}\|_{(V_L^\perp)'\to V_L^\perp}\le\beta_0^{-1},
 \qquad
 W[h]\ge\tfrac{\kappa_0}{2}\sigma .
\end{align}
Since $\det W[h(1)]=\rho$, \eqref{eq:quadratic} with $w=-\eta$ and
\eqref{eq:greens} express the residual at $x_0=h(1)-\eta$ through $\eta$ and
$a_L$, and the $b_z$ bound ($z=h(1),\eta$), \eqref{eq:approx}, and
\eqref{eq:loadbounds} show that it is small:
\begin{align}\label{eq:newtonresidual}
 G_L(x_0)=b_{h(1)}(\eta,\cdot)-\tfrac12b_\eta(\eta,\cdot)-(a_L,\cdot),
 \qquad
 \|G_L(x_0)\|_{(V_L^\perp)'}\le CL^{-(s+1)}.
\end{align}

\emph{Step 3: existence and uniqueness of the Galerkin root.}
By \eqref{eq:newtonlip}--\eqref{eq:newtonresidual}, the
Newton--Kantorovich parameter at $x_0$, the product of the
bounds $\beta_0^{-1}$ on $\|G_L'(x_0)^{-1}\|$, $CL^2$ on the Lipschitz
constant of $G_L'$, and $\beta_0^{-1}CL^{-(s+1)}$ on the first Newton
step, is $O(L^{1-s})\le\frac12$ for $L\ge L_1$.  So
\cite[Ch.~2]{deuflhard2011} gives a root $h_L^\ast$
with $\|h_L^\ast-x_0\|_{H^1}\le CL^{-(s+1)}\le\delta_1L^{-2}/4=R/4$.  Thus
$h_L^\ast\in\mathbb B_L$, and \eqref{eq:galerkinerror} follows from
\eqref{eq:newtoninverse} and
$\|h_L^\ast-h(1)\|_{H^1}\le\|h_L^\ast-x_0\|_{H^1}+\|\eta\|_{H^1}
\le CL^{-(s+1)}$.

For uniqueness, let $\tilde h_L^\ast\in\mathbb B_L$ be another root.
Steps~1 and~2 give
$0=G_L(\tilde h_L^\ast)-G_L(h_L^\ast)=G_L'\bigl((\tilde h_L^\ast+h_L^\ast)/2\bigr)(\tilde h_L^\ast-h_L^\ast)$, and
\eqref{eq:newtoninverse} at the midpoint, which lies in $\mathbb B_L$, gives $\tilde h_L^\ast=h_L^\ast$.

\emph{Step 4: quadratic convergence.}
Let $h\in\mathbb B_L$ and $e:=h-h_L^\ast$.  As $G_L(h_L^\ast)=0$,
Step~1 ($w=-e$) gives $G_L(h)=G_L'(h)e-(\det W[e],\cdot)$, so the new
error is
\[
 h^+-h_L^\ast=G_L'(h)^{-1}(\det W[e],\cdot),
 \qquad
 {\|h^+-h_L^\ast\|_{H^1}\le\beta_0^{-1}CL^2\|e\|_{H^1}^2
 =:CL^2\|e\|_{H^1}^2},
\]
by \eqref{eq:newtonlip}--\eqref{eq:newtoninverse}; this is
\eqref{eq:newtonquad}.  For $h\in\mathbb B_L$, Step~3 gives $\|e\|_{H^1}\le2R$,
so \eqref{eq:newtonquad} yields $\|h^+-h_L^\ast\|_{H^1}\le2C\delta_1\|e\|_{H^1}
\le\|e\|_{H^1}/4$ for $\delta_1$ small; hence
$\|h^+-x_0\|_{H^1}\le\|e\|_{H^1}/4+R/4<R$, and the Newton iterates stay in
$\mathbb B_L$ and converge quadratically to $h_L^\ast$.

Finally, the continuation ends in $\mathbb B_L$.  Since
$h^N-x_0=P_L(h^N-h(1))$, Theorems~\ref{thm:main} and~\ref{thm:trap} give
$\|h^N-x_0\|_{H^1}\le C(\tau^p+L^{-(s+1)})$, with $p=1$ for Euler and
$p=3$ for trapezoidal.  For large $L$, $CL^{-(s+1)}\le R/2$, and for $c$
small the conditions on $\tau$ give $\tau\le\tau_0$ and $C\tau^p\le R/2$.
\end{proof}

\section{Numerical experiments}\label{sec:experiments}

The implementation uses Python, NumPy, and SciPy; smooth integrals
use the quadrature rule $Q_L$ of \eqref{eq:QL} with
$n_x=2L+12$ and $n_\varphi=4L+13$, so that $E_L=4L+12$ exceeds the
requirement \eqref{eq:QLexact}; its nodes form the solver grid.
The final degree-$L$ state is denoted by $h_L$, either $h^N$ or
$h^{N,M}$ from Algorithm~\ref{alg:scheme}, as specified below.
All states satisfy $(h_L,\nu_k)=0$, $k=1,2,3$, fixing translation.
Unless stated otherwise, Newton correction stops when the
coefficient vectors $\mathbf h^m$ of $h^{N,m}$ satisfy
$\|\mathbf h^{m+1}-\mathbf h^m\|_2\le10^{-12}(1+\|\mathbf h^{m+1}\|_2)$.

All $H^1$ errors use the full norm of \eqref{eq:sobolev}.
For a benchmark whose exact support function $h^\ast=h(1)$ is known,
we evaluate
$e=h_L-h^\ast$ by
\begin{equation*}
 \|e\|_{H^1}
 \approx
 \bigl[Q_L(e^2+|\nabla e|^2)\bigr]^{1/2}.
\end{equation*}
Unless specified otherwise, $L^\infty$ residuals and minimum
eigenvalues of $W[h_L]$ are sampled on the solver grid;
positivity is the numerical convexity check.
Smooth strictly convex surfaces are plotted on a separate grid
using \eqref{eq:embedding},
$X_L(\nu)=h_L(\nu)\nu+\nabla h_L(\nu)$.
Here $\nabla h_L$ is evaluated by \eqref{eq:gradframe}; the
derivatives $\partial_\theta h_L$ and $\partial_\varphi h_L$ are
obtained by differentiating the spherical harmonic expansion of $h_L$
term by term, without finite differences.

\subsection{Ellipsoid benchmark}
\label{subsec:expA}\label{subsec:expB}

We use the ellipsoid $(a,b,c)=(1,1.5,2.5)$, whose support
function and density are
\begin{equation*}
 h^\ast(\nu)
 =\bigl(a^2\nu_1^2+b^2\nu_2^2+c^2\nu_3^2\bigr)^{1/2},
 \qquad
 \rho=\frac{(abc)^2}{(h^\ast)^4}.
\end{equation*}
The density satisfies \eqref{eq:closure} by symmetry.

\begin{figure}[htbp]
\centering
\begin{subfigure}[t]{0.48\linewidth}
\centering
\includegraphics[width=0.94\linewidth]{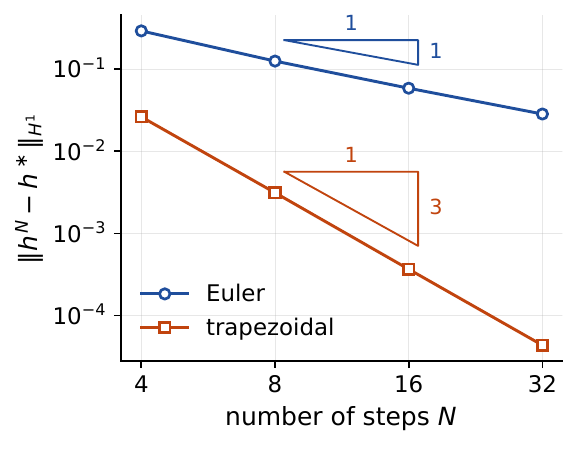}
\caption{Temporal convergence ($L=28$).}
\label{fig:convergence}
\end{subfigure}\hfill
\begin{subfigure}[t]{0.48\linewidth}
\centering
\includegraphics[width=0.94\linewidth]{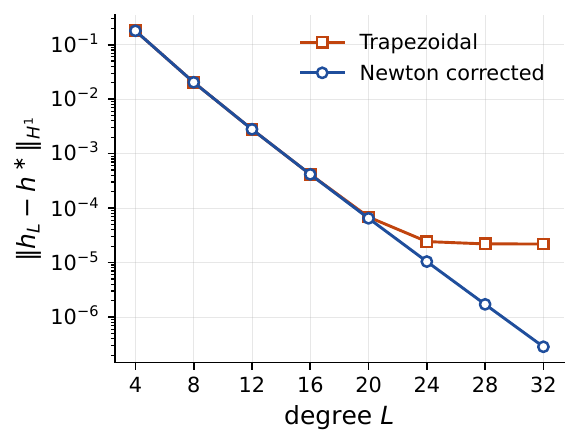}
\caption{Spectral convergence ($N=40$).}
\label{fig:spectral}
\end{subfigure}
\caption{Temporal and spectral $H^1$ convergence.}
\end{figure}

\paragraph{Temporal convergence.}
We fix $L=28$ and take $N=4,8,16,32$ steps without Newton
correction (Figure~\ref{fig:convergence}). We measure $e_N=\|h^N-h^\ast\|_{H^1}$ and compute
orders $\log_2(e_N/e_{2N})$. At this degree, spatial errors are
much smaller than temporal errors (Figure~\ref{fig:spectral}).
The orders are $1.22$, $1.10$, $1.04$ for Euler and $3.07$, $3.09$,
$3.08$ for the trapezoidal variant (Theorems~\ref{thm:main} and~\ref{thm:trap}).

\paragraph{Spectral convergence.}
We take $N=40$ trapezoidal steps with $L=4,8,\ldots,32$.
For $L=4$, quadrature uses $n_x=36$ and $n_\varphi=61$.
Figure~\ref{fig:spectral} plots $H^1$ errors against $h^\ast$:
orange before Newton correction and blue after correction.
Three updates suffice at each tested $L$.
For $L\le20$ the spatial error dominates, so the two curves
nearly coincide; for larger $L$ the orange curve levels off near
$2.19\times10^{-5}$, the time error for $N=40$
(Figure~\ref{fig:convergence} predicts $2.2\times10^{-5}$).
The blue curve, the error of the Galerkin solution $h_L^\ast$
of Proposition~\ref{prop:newton}, decreases to $2.85\times10^{-7}$ at
$L=32$.  Since $\rho$ is analytic here, the figure shows
spectral convergence.

\paragraph{Convexity and Newton initialization.}
At $L=20$, let $R_\infty(h)$ be the sampled maximum of
$|\det W[h]-\rho|$ and $\lambda_{\min}(W[h])$ the sampled minimum
eigenvalue. One Newton step from the sphere $h^{0,0}=1$ gives
\begin{align*}
 R_\infty(h^{0,0})&=1.29\times10^1,
 &
 R_\infty(h^{0,1})&=3.45\times10^1,
 &
 \lambda_{\min}(W[h^{0,1}])&=-0.346.
\end{align*}
The residual increases and convexity is lost. In contrast,
$N=20$ Euler continuation steps give
\begin{equation*}
 \min_{0\le n\le N}\lambda_{\min}(W[h^n])=0.396,
 \qquad
 R_\infty(h^N)=0.385.
\end{equation*}
Thus every continuation state passes the convexity check, and
$h^{N,0}=h^N$ supplies the initial state for Newton correction.

\paragraph{Compatibility.}
Adding $\varepsilon\nu_3$ to $\rho$, with
$\varepsilon\in\{0,0.02,0.05,0.10\}$, preserves positivity
but violates closure for $\varepsilon>0$:
\begin{equation*}
 \int_{\sph}\nu(\rho+\varepsilon\nu_3)\dA
 =\left(0,0,\frac{4\pi\varepsilon}{3}\right),
 \qquad
 \lambda=(0,0,-\varepsilon)
\end{equation*}
by \eqref{eq:lambdaformula}.
We use $L=20$, $N=20$ Euler steps, and Newton correction.
Table~\ref{tab:compatibility} reports the largest continuation multiplier and final
$L^\infty$ residuals against the perturbed input
($R_{\mathrm{in}}$) and the original density ($R_\rho$).

\begin{center}
\begin{minipage}{\linewidth}
\centering
\small
\setlength{\tabcolsep}{6pt}
\begin{tabular}{rcccc}
\toprule
$\varepsilon$
& predicted $\|\lambda\|_\infty$
& $\max_{0\le n<N}\|\lambda^n\|_\infty$
& $R_{\mathrm{in}}$
& $R_\rho$\\
\midrule
0 & 0 & $7.7\times10^{-16}$ & $3.265\times10^{-3}$ & $3.265\times10^{-3}$\\
0.02 & 0.02 & 0.02 & $2.097\times10^{-2}$ & $3.265\times10^{-3}$\\
0.05 & 0.05 & 0.05 & $5.094\times10^{-2}$ & $3.265\times10^{-3}$\\
0.10 & 0.10 & 0.10 & $1.009\times10^{-1}$ & $3.265\times10^{-3}$\\
\bottomrule
\end{tabular}
\captionsetup{hypcap=false}
\captionof{table}{Compatibility diagnostics for the perturbed ellipsoid density.}
\label{tab:compatibility}
\end{minipage}
\end{center}

The multiplier detects the incompatible component, which the
saddle formulation removes. Thus $R_{\mathrm{in}}$ grows with
$\varepsilon$, while $R_\rho$ stays at $3.265\times10^{-3}$:
all four inputs give the same numerical reconstruction to the digits shown.

\subsection{Analytic data without a closed-form support function}
\label{subsec:expS}

We prescribe the analytic density
\begin{equation}\label{eq:stardensity}
 \rho
 =
 \exp\left(
   \tfrac32\cos5\varphi\,\sin^5\theta
   +\tfrac3{10}P_2(\cos\theta)
 \right),
 \qquad
 \frac{\max\rho}{\min\rho}\approx20,
\end{equation}
where $P_2(x)=\tfrac12(3x^2-1)$.
Fivefold rotation and equatorial reflection ensure \eqref{eq:closure}.
We use $N=40$ trapezoidal steps followed by Newton correction.
With the corrected $L=32$ state $h_{32}$ as reference, the
table in Figure~\ref{fig:starbody} lists $\|h_L-h_{32}\|_{H^1}$, computed
exactly from the spherical harmonic coefficients by \eqref{eq:sobolev},
and $\lambda_{\min}(W[h_L])$ sampled on a $180\times360$ grid.

\begin{center}
\begin{minipage}[c]{0.61\linewidth}
\centering
\includegraphics[width=\linewidth]{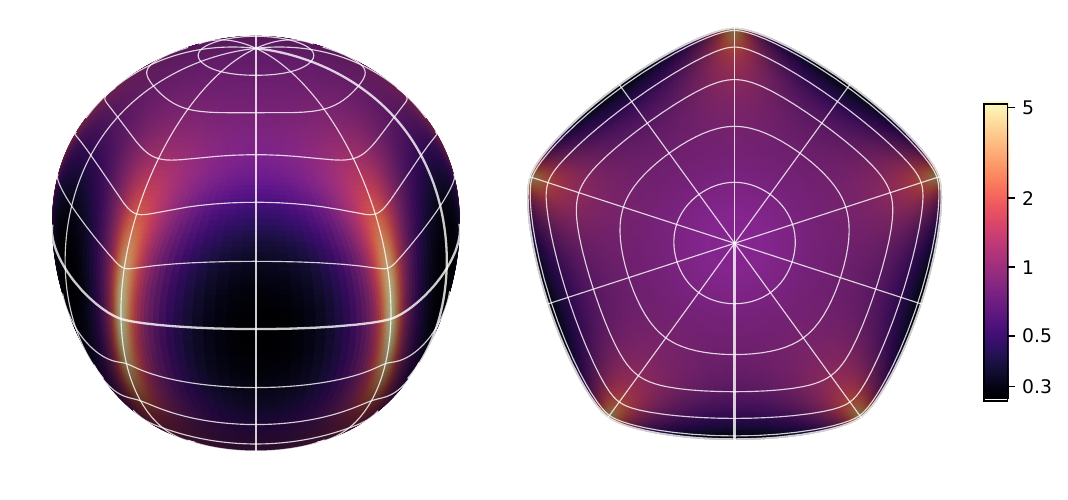}
\end{minipage}\hfill
\begin{minipage}[c]{0.37\linewidth}
\centering
\small
\setlength{\tabcolsep}{3pt}
\resizebox{\linewidth}{!}{%
\begin{tabular}{rcc}
\toprule
$L$ & $\lambda_{\min}(W[h_L])$
& $\|h_L-h_{32}\|_{H^1}$\\
\midrule
12 & 0.233   & $5.717\times10^{-3}$\\
16 & 0.180   & $5.014\times10^{-4}$\\
20 & 0.189   & $2.971\times10^{-5}$\\
\bottomrule
\end{tabular}%
}
\end{minipage}
\captionsetup{hypcap=false}
\captionof{figure}{Reconstruction of \eqref{eq:stardensity}
at $L=24$, colored by Gauss curvature $\mathcal K=1/\det W[h_L]$
on a logarithmic scale. Left: oblique view; right: view along the
fivefold symmetry axis.}
\label{fig:starbody}
\end{center}

\subsection{Reconstruction from discrete surface-area measures}
\label{subsec:expPolyhedra}

\begin{figure}[htbp]
\centering
\begin{subfigure}[t]{0.19\linewidth}
\centering
\includegraphics[width=\linewidth]{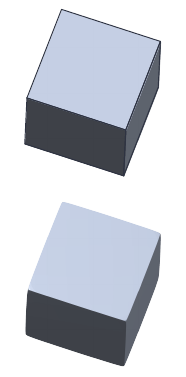}
\caption{}
\label{fig:polyhedraDirect:cube}
\end{subfigure}\hfill
\begin{subfigure}[t]{0.19\linewidth}
\centering
\includegraphics[width=\linewidth]{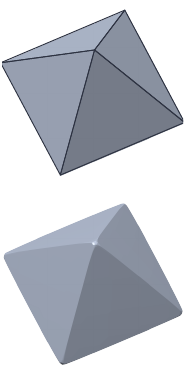}
\caption{}
\label{fig:polyhedraDirect:octahedron}
\end{subfigure}\hfill
\begin{subfigure}[t]{0.19\linewidth}
\centering
\includegraphics[width=\linewidth]{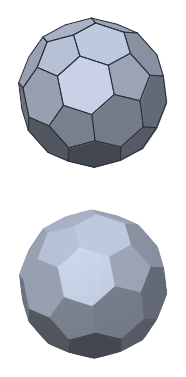}
\caption{}
\label{fig:polyhedraDirect:truncated_icosahedron}
\end{subfigure}\hfill
\begin{subfigure}[t]{0.19\linewidth}
\centering
\includegraphics[width=\linewidth]{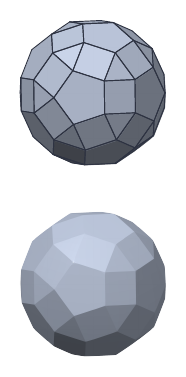}
\caption{}
\label{fig:polyhedraDirect:rhombicosidodecahedron}
\end{subfigure}\hfill
\begin{subfigure}[t]{0.19\linewidth}
\centering
\includegraphics[width=\linewidth]{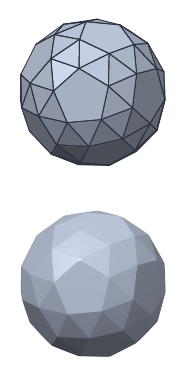}
\caption{}
\label{fig:polyhedraDirect:snub_dodecahedron}
\end{subfigure}
\caption{Exact polytopes (top) and finite half-space approximations of $K_L$
(bottom), $L=128$: (a)~cube; (b)~octahedron;
(c)~truncated icosahedron; (d)~rhombicosidodecahedron;
(e)~snub dodecahedron. Sampled Hausdorff-distance estimates are
$0.024$, $0.050$, $0.003$, $0.004$, and $0.007$.}
\label{fig:polyhedraDirect}
\end{figure}

A convex polytope $P$ is determined
up to translation by $\mu=\sum_f A_f\delta_{n_f}$, where $A_f$
and $n_f$ are its facet areas and outward unit normals,
and $\delta_{n_f}$ is the unit point mass at $n_f$
\cite[Chapter~8]{schneider2014}.
The five polytopes in Figure~\ref{fig:polyhedraDirect} are
rescaled to area $4\pi$ and generically rotated.
We use $\mu$ without mollification: the data term $Q_L(fq)$
in \eqref{eq:discretevelocity} is replaced by
$\sum_f A_fq(n_f)-(1,q)$.
Since $\rho$ is a sum of point masses and the support function
$h_P$ of $P$ is only Lipschitz, the hypotheses of
Theorem~\ref{thm:main} fail.

At $L=128$, we take $100$ Euler continuation steps and $12$
Newton updates. When a linear solve does not reach its tolerance,
we accept the current Krylov iterate. The discrete equations are
therefore not solved to convergence, and the final $W[h_L]$ is
indefinite.

A smooth $h$ is the support function of a possibly
degenerate convex body if and only if $W[h]\ge0$: its
homogeneous extension $\Phi(x)=|x|h(x/|x|)$ satisfies
$D^2\Phi|_{\sph}=\operatorname{diag}(W[h],0)$ (cf.\ the proof of
Lemma~\ref{lem:Winfty}), and convexity of $\Phi$ characterizes support
functions \cite[Section~1.7]{schneider2014}. Hence $h_L$ is not a
support function.  The planes $x\cdot\nu=h_L(\nu)$, $\nu\in\sph$,
still cut out a convex body,
\begin{equation}\label{eq:KL}
 K_L=\bigcap_{\nu\in\sph}
 \{x\in\mathbb R^3:x\cdot\nu\le h_L(\nu)\},
\end{equation}
whose support function $h_{K_L}$ is the largest
support function bounded by $h_L$.  In practice we use the normals
of a $400\times800$ grid.

Translation is fixed by $(h_L,\nu_k)=0$, $k=1,2,3$, and the Steiner
point of $P$ is at the origin. Since the
Hausdorff distance of convex bodies is the maximum norm of the
difference of their support functions \cite[Section~1.8]{schneider2014},
Figure~\ref{fig:polyhedraDirect} reports $\max|h_{K_L}-h_P|$ over a
separate $320\times640$ grid of directions.

\subsection{Umbilic locations and indices}\label{subsec:expU}

We test recovery of umbilic locations and indices using $N=20$
Euler steps and Newton correction. In the spherical orthonormal frame
$(e_\theta,e_\varphi)$ of \eqref{eq:gradframe}, set
\begin{equation*}
 q=(W_{11}-W_{22})+2\mathrm{i}W_{12},
 \qquad |q|=|r_1-r_2|,
\end{equation*}
where $r_1,r_2$ are the eigenvalues of $W[h_L]$ and $\mathrm i^2=-1$.
Away from the poles, where the frame is singular,
umbilics are the zeros of $q$. If $\alpha$ is the angle of a principal
direction, then $q=\pm|q|e^{2\mathrm i\alpha}$, so the phase of $q$
winds twice as fast as the principal directions.
Around each umbilic, we sample a small loop $\Gamma$
on which $q\ne0$, unwrap its phase, and set
$w=\frac{1}{2\pi}\oint_\Gamma\mathrm{d}\arg q$.
Since the frame itself turns once around each pole,
\begin{equation*}
 \operatorname{ind}=\tfrac12w,
 \qquad
 \operatorname{ind}_N=1+\tfrac12w,
 \qquad
 \operatorname{ind}_S=1-\tfrac12w.
\end{equation*}
Nonpolar loops are positively oriented; polar latitude loops
follow increasing $\varphi$.
Sampling is refined until $w$ stabilizes.

\paragraph{Spheroid.}
For $(a,b,c)=(1,1,1.35)$, the support function and
the meridional--azimuthal radius difference are
\begin{equation*}
 h^\ast=(a^2\sin^2\theta+c^2\cos^2\theta)^{1/2},
 \qquad
 r_\theta-r_\varphi
 =\frac{a^2(c^2-a^2)\sin^2\theta}{(h^\ast)^3}.
\end{equation*}
Only the poles are umbilic; at $L=24$, both computed indices
are $+1$ ($w=0$). The maximum absolute error of
the computed $W_{11}-W_{22}$ against this formula, over $500$
equispaced polar angles in $[0.02,\pi-0.02]$ at $\varphi=0$, is
$3.9\times10^{-10}$.

\paragraph{Triaxial ellipsoid.}
For $(a,b,c)=(1,1.2,1.5)$, the four umbilics lie in
the $xz$-plane, with exact normals $(\pm0.742,0,\pm0.671)$ to three
decimals, signs independent.
Without using this information, a $240\times480$
full-sphere search at $L=24$ detects four candidates as local minima
of $|q|$ below $8\%$ of its sampled maximum
(Figure~\ref{fig:umbilicGap}).
To measure convergence, we then compute $q$ from the
spherical harmonic expansion of $h_L$ and find its zeros
by Newton's method with the exact normals as initial
guesses.
Table~\ref{tab:umbilicLocations} reports the largest angle between
computed and exact normals over the four umbilics;
it decreases by more than a factor of $10$ whenever $L$
increases by $4$, consistent with spectral convergence.
All four indices are $+1/2$ ($w=1$), summing to the Euler
characteristic $\chi(\sph)=2$.

\begin{center}
\begin{minipage}{\linewidth}
\centering
\small
\setlength{\tabcolsep}{5pt}
\begin{tabular}{lccccc}
\toprule
$L$ & 8 & 12 & 16 & 20 & 24\\
Max. angular error (rad)
& $4.2\times10^{-4}$
& $3.8\times10^{-5}$
& $2.5\times10^{-6}$
& $1.4\times10^{-7}$
& $7.4\times10^{-9}$\\
\bottomrule
\end{tabular}
\captionsetup{hypcap=false}
\captionof{table}{Convergence of umbilic locations on the triaxial ellipsoid.}
\label{tab:umbilicLocations}
\end{minipage}
\end{center}

\begin{center}
\begin{minipage}{\linewidth}
\captionsetup{type=figure,skip=4pt}
\centering
\begin{subfigure}[b]{0.48\linewidth}
\centering
\includegraphics[width=0.64\linewidth]{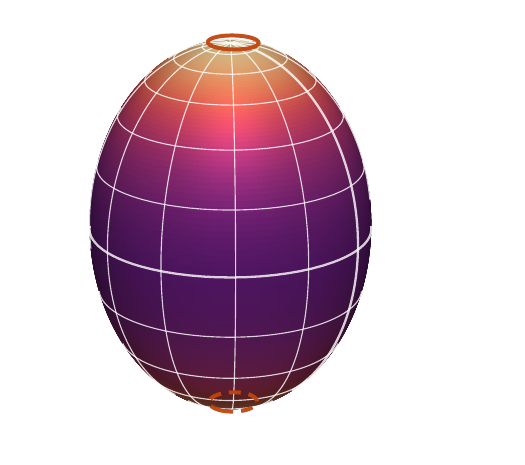}
\caption{Spheroid.}
\end{subfigure}\hfill
\begin{subfigure}[b]{0.48\linewidth}
\centering
\includegraphics[width=0.64\linewidth]{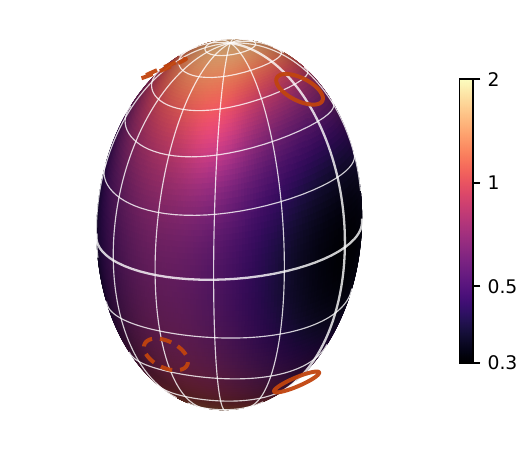}
\caption{Triaxial ellipsoid.}
\end{subfigure}
\caption{Reconstructions at $L=24$, colored by Gauss curvature
$\mathcal{K}$ on a shared logarithmic scale.  Circles mark exact umbilics;
polar markers are slightly displaced.}
\label{fig:umbilics}
\end{minipage}
\end{center}

Figure~\ref{fig:umbilics} shows the $L=24$ reconstructions.
For Figure~\ref{fig:umbilicPhase}, we evaluate $q$ directly
on the geodesic circle $\Gamma$ of radius $0.08$ centered at
$(\theta,\varphi)\approx(0.8355,0)$, using $720$ equal intervals.
The sampled minimum of $|q|$ is $0.067$.
Appendix~\ref{app:phase} gives $|q|\ge0.055$ on $\Gamma$ and
bounds each phase increment by $0.45<\pi$. Thus phase unwrapping is
unambiguous; the total increase is $2\pi$, giving index $+1/2$.

\begin{center}
\begin{minipage}{\linewidth}
\captionsetup{type=figure,skip=4pt}
\centering
\begin{subfigure}[b]{0.48\linewidth}
\centering
\includegraphics[width=0.76\linewidth]{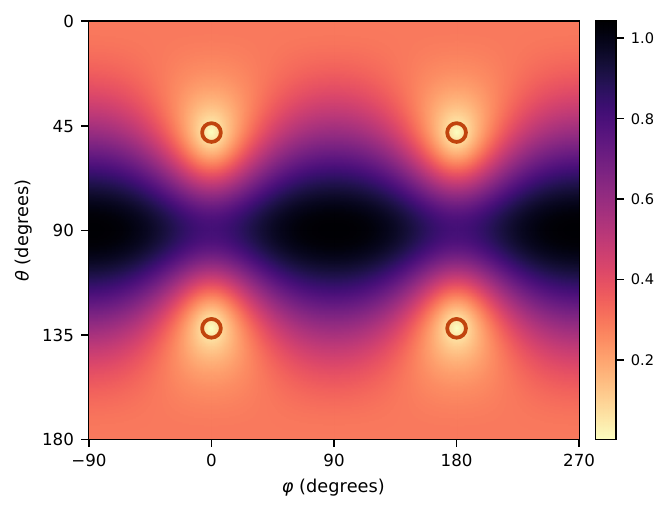}
\caption{Principal-radius difference.}
\label{fig:umbilicGap}
\end{subfigure}\hfill
\begin{subfigure}[b]{0.48\linewidth}
\centering
\includegraphics[width=0.76\linewidth]{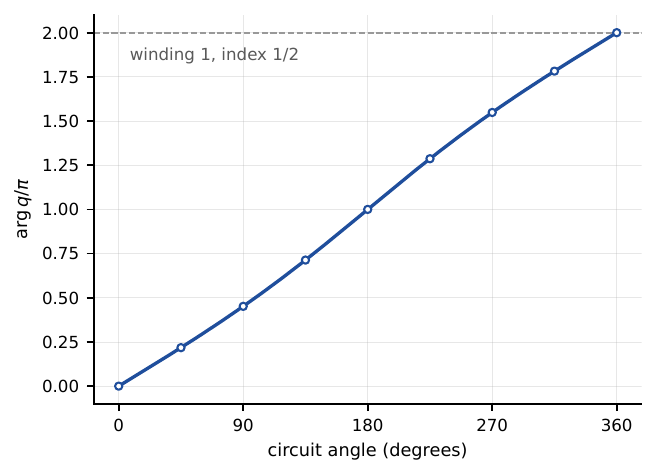}
\caption{Unwrapped phase.}
\label{fig:umbilicPhase}
\end{subfigure}
\caption{Umbilic diagnostics for the triaxial ellipsoid at $L=24$.
(a)~$|r_1-r_2|$ versus normal direction;
circles mark exact umbilics.
(b)~Unwrapped phase of $q$ around one umbilic.}
\label{fig:umbilicDiag}
\end{minipage}
\end{center}

\section{Conclusion}

\nobreak
We have analyzed a spherical harmonic
continuation from the unit sphere that
reduces the Minkowski problem to linear elliptic equations
in divergence form.  For positive, compatible $\rho\in H^{s+2}(\sph)$ with
real $s\ge4$, the Euler and trapezoidal states stay strictly convex and
have $H^1$ errors $O(\tau+L^{-(s+1)})$ and $O(\tau^{3}+L^{-(s+1)})$,
including quadrature, with independent thresholds for $\tau$ and
$L$.  The proof uses a priori bounds for each iterate's perturbed equation,
which restore ellipticity at every step, and the exactly quadratic
determinant; the multiplier measures
the incompatibility of the data.  Under $\tau\le cL^{-2}$ (Euler) or
$\tau\le cL^{-2/3}$ (trapezoidal), the final state lies in a ball
where Newton's method converges quadratically to the Galerkin
solution; started directly from the sphere,
the first undamped step left the convex cone
(Section~\ref{subsec:expA}).  Experiments illustrate the rates and the compatibility
diagnostic, reconstruct convex polyhedra from facet normals and areas
beyond the theory, and recover umbilics.

\bibliographystyle{siam}
\bibliography{refs}

\appendix
\section*{Appendix}
\addcontentsline{toc}{section}{Appendix}
\stepcounter{section}

\subsection{Proofs for Section~\ref{sec:formulation}}\label{app:two}

\begin{proof}[Proof of Lemma~\ref{lem:divfree}]
Since $\nabla\sigma=0$, differentiating \eqref{eq:Wdef} gives
\begin{align*}
\nabla_iW[h]_{jk}=\nabla_i\nabla_j\nabla_kh+(\nabla_ih)\,\sigma_{jk},
\end{align*}
and antisymmetrizing in $i$ and $j$,
\begin{align*}
\nabla_iW[h]_{jk}-\nabla_jW[h]_{ik}
=[\nabla_i,\nabla_j]\nabla_kh
+\sigma_{jk}\nabla_ih-\sigma_{ik}\nabla_jh .
\end{align*}
Since $\sph$ has constant curvature $1$, the Ricci identity reads
\begin{align*}
[\nabla_i,\nabla_j]\nabla_kh=\sigma_{ik}\nabla_jh-\sigma_{jk}\nabla_ih .
\end{align*}
This cancels the last two terms above, so
$\nabla_iW[h]_{jk}=\nabla_jW[h]_{ik}$, that is, $W[h]$ is a
Codazzi tensor.
To deduce \eqref{eq:divfree}, contract this identity
with $\sigma^{ik}$ and use $\nabla\sigma=0$:
\begin{align*}
\nabla^kW[h]_{jk}=\nabla_j\tr W[h],
\end{align*}
that is, the divergence of $W[h]$ equals the gradient of its trace.
Hence, by \eqref{eq:cofactor} and the symmetry of $W[h]$,
\begin{align*}
\nabla_iC[h]^{ij}
=\nabla_i\bigl((\tr W[h])\,\sigma^{ij}-W[h]^{ij}\bigr)
=\nabla^j\tr W[h]-\nabla^j\tr W[h]
=0 . \tag*{\qedhere}
\end{align*}
\end{proof}

\begin{proof}[Verification of the symmetry of $T$ in \eqref{eq:trilinear}]
Since the
transpositions of $(g,w)$ and of $(f,w)$ generate all permutations of the
three arguments, it suffices to verify these two. Symmetry in $(g,w)$
holds pointwise: it is the symmetry \eqref{eq:mixeddisc} of the mixed
discriminant. Symmetry in $(f,w)$ is the formal
self-adjointness of $\mathcal{L}_g$: for $g\in C^3(\sph)$, by \eqref{eq:Lop} and
\eqref{eq:greens},
\begin{align*}
T(f,g,w)
=\int_{\sph}f\,\mathcal{L}_gw\dA
=-b_g(w,f)
=-b_g(f,w)
=\int_{\sph}w\,\mathcal{L}_gf\dA
=T(w,g,f),
\end{align*}
where the second and fourth equalities are \eqref{eq:greens} and the
third is the symmetry of $b_g$.  For $g\in C^2(\sph)$ we
approximate: since
\begin{align*}
|T(f,g,w)|\le C\|f\|_{C^0}\|g\|_{C^2}\|w\|_{C^2},
\end{align*}
$T$ is continuous on $C^2(\sph)^3$, and $C^\infty(\sph)$ is dense in
$C^2(\sph)$, so the symmetry passes to the limit.
Theorem~\ref{thm:kernel} uses this case, with $g=v\in C^{2,\alpha}(\sph)$.
\end{proof}

\begin{proof}[Proof of Lemma~\ref{lem:pointwise}]
An orthogonal change of frame, which leaves $D$ and $\det$
invariant, gives $A=\mathrm{diag}(\mu_1,\mu_2)$ with $\mu_1,\mu_2>0$.
Then $D(A,B)=\tfrac12(\mu_1B_{22}+\mu_2B_{11})$
and
\begin{align*}
D(A,B)^2-\det A\det B
=\tfrac14\bigl(\mu_1B_{22}-\mu_2B_{11}\bigr)^2
+\mu_1\mu_2B_{12}^2\;\ge\;0 .
\end{align*}
If $D(A,B)=0$, this and $\det A>0$ give $\det B\le0$. If
moreover $\det B=0$, the display forces $\mu_1B_{22}=\mu_2B_{11}$ and
$B_{12}=0$; with $\mu_1B_{22}+\mu_2B_{11}=2D(A,B)=0$ this gives $B=0$.
\end{proof}

\begin{proof}[Proof of Theorem~\ref{thm:kernel}]
Since $W[\nu_k]=0$, the inclusion $\mathcal{N}\subseteq\ker\mathcal{L}_h$
holds for every $h$. Conversely, let $v\in H^2(\sph)$ satisfy $\mathcal{L}_hv=0$. Since
$h\in C^3(\sph)$ and $W[h]>0$, the coefficient $C[h]$ is of class $C^1$ and uniformly
elliptic, so elliptic regularity gives $v\in C^{2,\alpha}(\sph)$ for every
$\alpha<1$. By
\eqref{eq:Lop}, the symmetry of $T$ in
\eqref{eq:trilinear}, whose $C^2$ case covers
$v\in C^{2,\alpha}$, and \eqref{eq:cofactorids},
\begin{align*}
0=\int_{\sph} v\,\mathcal{L}_hv\dA
=T(v,h,v)=T(h,v,v)
=2\int_{\sph}h\,\det W[v]\dA .
\end{align*}
At every point, $D(W[h],W[v])=\tfrac12\mathcal{L}_hv=0$ and $W[h]>0$, so
Lemma~\ref{lem:pointwise} gives $\det W[v]\le0$. Since $h>0$, the
continuous function $h\det W[v]\le0$ has zero integral, so
$\det W[v]\equiv0$, and Lemma~\ref{lem:pointwise}
then gives $W[v]\equiv0$. Taking the trace, $\Delta v+2v=0$, so $v$ is a
spherical harmonic of degree one, that is, $v\in\mathcal{N}$.
\end{proof}

\begin{proof}[Proof of Proposition~\ref{prop:pathreg}]
Step 1: existence. Each $\rho_t\in C^{k,\alpha}(\sph)$ is
admissible by Lemma~\ref{lem:path}, so by Theorem~\ref{thm:classical}
the solution $h(t)$ of \eqref{eq:pathproblem} exists, lies in
$C^{k+2,\alpha}(\sph)$, and is positive.

Step 2: the linearization is an isomorphism. Set
\begin{align*}
X&=\{h\in C^{k+2,\alpha}(\sph):(h,\nu_j)=0,\;j=1,2,3\},\\
Y_0&=\{f\in C^{k,\alpha}(\sph):(f,\nu_j)=0,\;j=1,2,3\}.
\end{align*}
Fix $t$. Since $W[h(t)]>0$, also $C[h(t)]>0$, so
$\mathcal{L}_{h(t)}$ is uniformly elliptic with
$C^{k,\alpha}$ coefficients.  Elliptic Schauder theory
\cite{gilbargtrudinger2001} makes it Fredholm of index zero
from $C^{k+2,\alpha}(\sph)$ to $C^{k,\alpha}(\sph)$.  By
\eqref{eq:greens}, its range lies in the $L^2$-orthogonal complement of
its kernel; both have codimension $\dim\ker\mathcal{L}_{h(t)}$, so they
coincide.  By Step~1,
Theorem~\ref{thm:kernel} applies to $h(t)$; hence
\begin{align*}
\ker\mathcal{L}_{h(t)}=\mathcal{N},
\qquad
\operatorname{ran}\mathcal{L}_{h(t)}=Y_0 .
\end{align*}
The kernels $X$ and $Y_0$ of the bounded projection
$P_{\mathcal N}$ are closed, and $C^{k+2,\alpha}(\sph)=X\oplus\mathcal{N}$;
so $\mathcal{L}_{h(t)}:X\to Y_0$ is a bounded bijection between Banach
spaces, hence an isomorphism.

Step 3: smooth dependence on $t$.  With $\rho_t$ given by
\eqref{eq:path} for all real $t$, define
\begin{align*}
\mathcal F:\R\times X\to Y_0,
\qquad
\mathcal F(t,h)=\det W[h]-\rho_t ,
\end{align*}
with values in $Y_0$ by Lemma~\ref{lem:closureid},
\eqref{eq:closure}, and $\int_{\sph}\nu\dA=0$. By \eqref{eq:cofactorids},
\begin{align*}
\mathcal F(t,h+u)-\mathcal F(t,h)=\mathcal{L}_{h}u+\det W[u],
\qquad u\in X,
\end{align*}
without remainder. Since
$\|fg\|_{C^{k,\alpha}(\sph)}\le C\|f\|_{C^{k,\alpha}(\sph)}\|g\|_{C^{k,\alpha}(\sph)}$
and the entries of $W[u]$ lie in $C^{k,\alpha}(\sph)$ for $u\in X$,
$u\mapsto\det W[u]$ is a bounded quadratic form from $X$ to
$C^{k,\alpha}(\sph)$, so
$\mathcal F$ is a continuous quadratic
polynomial in $h$ (and affine in $t$), hence jointly $C^\infty$ in $(t,h)$ with
$D_h\mathcal F(t,h)=\mathcal{L}_{h}$.

Fix $t_0\in[0,1]$.
Since $\mathcal F(t_0,h(t_0))=0$ and $\mathcal{L}_{h(t_0)}:X\to Y_0$ is
an isomorphism by Step~2, the implicit function theorem gives a
$C^\infty$ map $g$ on a neighborhood~$I$ of $t_0$ with
$\mathcal F(t,g(t))=0$ and $g(t_0)=h(t_0)$; shrinking $I$, also
$W[g(t)]>0$ by continuity.  For $t\in I\cap[0,1]$, $g(t)$ thus
solves \eqref{eq:pathproblem}, so $g(t)=h(t)$ by the uniqueness in
Theorem~\ref{thm:classical}.  Covering $[0,1]$ by
finitely many such $I$ gives
$h\in C^\infty\bigl([0,1];C^{k+2,\alpha}(\sph)\bigr)$.

Step 4: uniform bounds.  By Step~3, $W[h(t)]$ and $h(t)$
depend continuously on $t$ in $C(\sph)$.  The eigenvalues of $W[h(t)]$
and the values of $h(t)$ are therefore continuous positive functions on
the compact set $[0,1]\times\sph$, which gives $\kappa_0$, $\Lambda_0$,
and $h_{\min}$ in \eqref{eq:uniformbounds}; they depend only on $\rho$.
\end{proof}

\begin{proof}[Proof of Proposition~\ref{prop:multiplier}]
(i) Take $q=\nu_j$ in the first equation.  Since $W[\nu_j]=0$,
\eqref{eq:greens} and the symmetry of $b_h$ give
$b_h(w,\nu_j)=-(\mathcal{L}_h\nu_j,w)=0$ for every $w\in H^1(\sph)$,
and $(\nu_k,\nu_j)=\tfrac{4\pi}{3}\delta_{kj}$ by
\eqref{eq:kernelspace}, so $\lambda_j\,\tfrac{4\pi}{3}=-(f,\nu_j)$.
This is \eqref{eq:lambdaformula}, and $\lambda\cdot\nu=-P_{\mathcal{N}}f$
by \eqref{eq:PNdef}.  With this $\lambda$, the first equation reads
$b_h(v,q)=-(f-P_{\mathcal{N}}f,q)$.  Both sides vanish for
$q\in\mathcal{N}$, so it suffices to take $q\in\Hperp$, and
Proposition~\ref{prop:infsup} gives exactly one such $v\in\Hperp$.

(ii) By (i) with $h=\bar h(t)$ and $f=\partial_t\bar\rho_t$,
$\lambda(t)\cdot\nu=-P_{\mathcal{N}}\partial_t\bar\rho_t$.  Since
$\bar h(t)$ is differentiable in $t$ with values in $H^2(\sph)$, the
product rule applies to $\det W[\bar h(t)]$, and by
\eqref{eq:cofactorids} its derivative is the linear term
$\mathcal{L}_{\bar h(t)}\partial_t\bar h(t)$.  With (i), this gives
\begin{align*}
\partial_t\det W[\bar h(t)]
=\mathcal{L}_{\bar h(t)}\partial_t\bar h(t)
=\partial_t\bar\rho_t+\lambda(t)\cdot\nu
=\partial_t\bigl(\bar\rho_t-P_{\mathcal{N}}\bar\rho_t\bigr).
\end{align*}
The two sides of \eqref{eq:wrongtarget} therefore have equal
$t$-derivatives, and they agree at $t=0$ because
$\det W[\bar h(0)]=1=\bar\rho_0-P_{\mathcal{N}}\bar\rho_0$; this proves
\eqref{eq:wrongtarget}. For the linear path $\bar\rho_t=\rho_t$, Lemma~\ref{lem:path} gives
$\partial_t\rho_t=\rho-1\perp\mathcal{N}$, so
$\lambda\equiv0$ by \eqref{eq:lambdaformula}, and
$P_{\mathcal{N}}\rho_t=0$ reduces \eqref{eq:wrongtarget} to
$\det W[\bar h(t)]=\rho_t$.
\end{proof}

\begin{proof}[Proof of Proposition~\ref{prop:infsup}]
Step 1: the operator is Fredholm of index zero. Define
$\mathcal A,\mathcal J:\Hperp\to(\Hperp)'$ by
$\langle\mathcal Av,q\rangle=b_h(v,q)$ and
$\langle\mathcal Jv,q\rangle=(v,q)$; $\mathcal A$ is bounded, and
$\mathcal J$ is compact by the Rellich embedding.
As $W[h]>0$ is continuous, $h\in\mathcal W(\kappa,\Lambda)$
for some $0<\kappa\le\Lambda$, and Lemma~\ref{lem:garding} gives
$\langle(\mathcal A+c_\ast\mathcal J)v,v\rangle\ge\kappa\|v\|_{H^1}^2$
with $c_\ast=\kappa+2\Lambda$.  So $\mathcal A+c_\ast \mathcal J$ is
invertible by Lax--Milgram, and $\mathcal A$, a compact
perturbation of it, is Fredholm of index zero.

Step 2: $\mathcal A$ is injective. Let $b_h(v,q)=0$ for all
$q\in \Hperp$. Since $b_h(v,\nu_k)=-(\mathcal{L}_h\nu_k,v)=0$
by \eqref{eq:greens}, in fact $b_h(v,q)=0$ for all
$q\in H^1(\sph)$, that is, $\mathcal{L}_hv=0$ weakly. Elliptic
regularity gives $v\in H^2(\sph)$, and Theorem~\ref{thm:kernel} gives
$v\in\mathcal{N}\cap \Hperp=\{0\}$.

Step 3: the inf-sup constant. For $g\in C^3(\sph)$ let
$\beta(g)$ be the left-hand side of \eqref{eq:infsup} with $b_g$ for $b_h$,
the best constant there.  By Steps~1--2, $\mathcal A$ is an
isomorphism, so for every $v\in \Hperp$,
\begin{align*}
\sup_{q\in \Hperp\setminus\{0\}}\frac{b_h(v,q)}{\|q\|_{H^1}}
=\|\mathcal Av\|_{(\Hperp)'}
\;\ge\;\frac{\|v\|_{H^1}}{\|\mathcal A^{-1}\|}\,,
\end{align*}
and hence $\beta(h)\ge1/\|\mathcal A^{-1}\|>0$.

Step 4: a perturbation bound.  Let $\tilde h\in C^3(\sph)$ and
$E=W[\tilde h]-W[h]$ with $\|E\|_{L^\infty}\le\varepsilon$.
Since $b_h$ depends on $h$ only through $W[h]$, and
linearly (through $C(W[h])$ and $\tr W[h]$ in \eqref{eq:greens}),
$b_{\tilde h}-b_h$ is the integral in \eqref{eq:greens} with $W[h]$
replaced by $E$.  Pointwise, $\|C(E)\|\le\|E\|$ since $C$ interchanges
the eigenvalues, and $|\tr E|\le2\|E\|$; by Cauchy--Schwarz in $L^2$ and
in $\R^2$,
\begin{align*}
|b_{\tilde h}(v,q)-b_h(v,q)|
\le\varepsilon\bigl(\|\nabla v\|_{L^2}\|\nabla q\|_{L^2}
+2\|v\|_{L^2}\|q\|_{L^2}\bigr)
\le2\varepsilon\,\|v\|_{H^1}\|q\|_{H^1}.
\end{align*}
So each supremum in \eqref{eq:infsup} drops by at most
$2\varepsilon$, and $\beta(\tilde h)\ge\beta(h)-2\varepsilon$.

Step 5: proof of (i) and (ii).  For $s,t\in[0,1]$, Step~4
applied in both directions gives
$|\beta(h(s))-\beta(h(t))|\le2\|W[h(s)]-W[h(t)]\|_{L^\infty}$, and
the right side tends to $0$ as $s\to t$ by Proposition~\ref{prop:pathreg}.
So $t\mapsto\beta(h(t))$ is continuous and positive on $[0,1]$, and its
minimum $\beta>0$ proves (i).  Step~4 with $h=h(t)$ and
$\beta(h(t))\ge\beta$ gives (ii).
\end{proof}

\subsection{Proofs for Section~\ref{sec:method}}\label{app:three}

\begin{proof}[Proof of Lemma~\ref{lem:Winfty}]
Step 1: $W[q]$ is the tangential Hessian of $\Phi$.
Let $\Phi(x)=|x|\,q(x/|x|)$ be the degree-one homogeneous extension of
$q$ to $\R^3\setminus\{0\}$.  Since $\Phi(r\nu)=r\,q(\nu)$,
\(
 \nabla_{\R^3}\Phi(r\nu)=\nabla q(\nu)+q(\nu)\,\nu ,
\)
which is independent of $r>0$; hence $D^2\Phi\,\nu=0$.  For an orthonormal
tangent frame $(e_1,e_2)$ at $\nu$, let $D_j$ denote the derivative along
$e_j$.  Then $D_j\nu=e_j$ and $D_jX=\nabla_jX-(X\cdot e_j)\,\nu$ for
tangent fields $X$, so
\begin{align*}
 D^2\Phi\,e_j
 =D_j\bigl(\nabla q+q\,\nu\bigr)
 =(\nabla_j\nabla^{i}q)\,e_i+q\,e_j
 =W[q]_j^{\;i}\,e_i .
\end{align*}
Hence, in Cartesian coordinates $x_1,x_2,x_3$ of $\R^3$
($\partial_a=\partial/\partial x_a$, $a=1,2,3$), with the orthogonal matrix
$R=(e_1\;e_2\;\nu)$, the Hessian $D^2\Phi=(\partial_a\partial_b\Phi)$ and
the gradient $\nabla_{\R^3}\Phi=(\partial_a\Phi)$ satisfy on $\sph$
\begin{align*}
 D^2\Phi=R\begin{pmatrix} W[q] & 0\\ 0 & 0\end{pmatrix}R^{-1},
 \qquad
 \nabla_{\R^3}\Phi=\nabla q+q\,\nu .
\end{align*}

Step 2: the degree count.
Each component $q_l\in\mathcal H_l$ of $q$ is the restriction of a
homogeneous harmonic polynomial $P_l$ of degree $l$, so
$q_l(x/|x|)=|x|^{-l}P_l(x)$ and
\begin{align*}
 \Phi=\sum_{l\le L}|x|^{\,1-l}P_l(x).
\end{align*}
Since $\partial_a(|x|^{k}P)=k\,|x|^{k-2}x_aP+|x|^{k}\partial_aP$, each
derivative preserves this form and raises the degree of the polynomial
factors by at most one.  On $\sph$, where $|x|=1$, a polynomial of
degree at most $d$ restricts to an element of
$\bigoplus_{m\le d}\mathcal{H}_m$.  Hence
\begin{align*}
 \partial_a\Phi|_{\sph}\in\bigoplus_{m\le L+1}\mathcal{H}_m,
 \qquad
 \partial_a\partial_b\Phi|_{\sph}\in\bigoplus_{m\le L+2}\mathcal{H}_m .
\end{align*}

Step 3: the pointwise bounds.
Apply \eqref{eq:nikolskii} to each Cartesian entry of $W[q]$ and sum
the squares, then use $\|W[q]\|_{L^2}\le C\|q\|_{H^2}$ and the inverse
inequality \eqref{eq:inverse}:
\begin{align*}
 \|W[q]\|_{L^\infty}
 \le\frac{L+3}{\sqrt{4\pi}}\,\|W[q]\|_{L^2}
 \le C\,L\,\|q\|_{H^2}
 \le C\,L^{2}\,\|q\|_{H^1}.
\end{align*}
For the second bound, apply \eqref{eq:nikolskii} to $q$ and to the
entries of $\nabla q+q\,\nu$, which have degree at most $L+1$ and
satisfy $\|\nabla q+q\,\nu\|_{L^2}=\|q\|_{H^1}$ since
$\nabla q\perp\nu$; as $|\nabla q|\le|\nabla q+q\,\nu|$, this gives
\(
 \|q\|_{L^\infty}+\|\nabla q\|_{L^\infty}\le C\,L\,\|q\|_{H^1}.
\)
\end{proof}

\begin{proof}[Proof of Proposition~\ref{prop:quadratureload}]
Split $\rho=P_L\rho+\rho_L^\perp$.  For $q\in V_L$ the product
$(P_L\rho{-1})\,q$ has degree at most $2L\le E_L$
by \eqref{eq:QLexact} and is integrated exactly by
$Q_L$, and $(\rho_L^\perp,q)=0$ since $\rho_L^\perp\perp V_L$;
so
\begin{align*}
 (a_L,q)=Q_L(\rho_L^\perp q).
\end{align*}
Since the weights of $Q_L$ are positive and $Q_L$ is exact for $1$
and $q^2$, the Cauchy--Schwarz inequality gives
\begin{align*}
 Q_L(|q|)\le Q_L(1)^{1/2}Q_L(q^2)^{1/2}=\sqrt{4\pi}\,\|q\|_{L^2},
\end{align*}
and hence
\begin{align*}
 |(a_L,q)|
 \le\|\rho_L^\perp\|_{L^\infty}\,Q_L(|q|)
 \le\sqrt{4\pi}\,\|\rho_L^\perp\|_{L^\infty}\|q\|_{L^2}.
\end{align*}
For the tail, write $\rho_L^\perp=\sum_{\ell>L}f_\ell$ with
$f_\ell\in\mathcal H_\ell$.  The addition theorem gives
$\|f_\ell\|_{L^\infty}\le\bigl((2\ell+1)/4\pi\bigr)^{1/2}\|f_\ell\|_{L^2}$,
so the Cauchy--Schwarz inequality and \eqref{eq:sobolev} yield
\begin{align*}
 \|\rho_L^\perp\|_{L^\infty}
 \le C\Bigl(\sum_{\ell>L}(2\ell+1)\bigl(1+\ell(\ell+1)\bigr)^{-(s+2)}\Bigr)^{1/2}
 \|\rho\|_{H^{s+2}}
 \le CL^{-(s+1)}\|\rho\|_{H^{s+2}}.
\end{align*}
As $a_L\in V_L$, the supremum over $q\in V_L$ gives the $L^2$ bound.
The inverse inequality \eqref{eq:inverse} with $p'=0$ then gives the
$H^p$ bounds, and $\|q\|_{L^2}\le\|q\|_{H^1}$ gives the dual bound.
Finally, $H^{3+2\alpha}(\sph)\hookrightarrow C^{2,\alpha}(\sph)$ and
the $H^p$ bound with $p=3+2\alpha\le s$ give the last bound.
\end{proof}

\subsection{Proofs for Section~\ref{sec:analysis}}\label{app:four}

\begin{lemma}[product and commutator estimates at real order]
\label{lem:realorder}
Let $p\ge2$ be real, $0<\alpha<1$, and $J^{p}:=(I-\Delta)^{p/2}$.
For functions on $\sph$, whenever the right-hand sides are finite,
\begin{align}
\|fg\|_{H^{p}}&\le C\bigl(\|f\|_{H^{p}}\|g\|_{C^{\alpha}}
+\|f\|_{C^{\alpha}}\|g\|_{H^{p}}\bigr),
\label{eq:moser}\\
\|J^{p}(au)-aJ^{p}u\|_{L^2}&\le C\bigl(\|\nabla a\|_{L^\infty}\|u\|_{H^{p-1}}
+\|a\|_{H^{p}}\|u\|_{L^\infty}\bigr).
\label{eq:katoponce}
\end{align}
Moreover, Lemma~\ref{lem:lineartame} and the estimate
\eqref{eq:localtame} in its proof hold with any real $p\ge2$ in place
of $s$, with a constant depending on $p$, on the coefficient bounds,
and, for Lemma~\ref{lem:lineartame}, on $\beta$.
\end{lemma}

\begin{proof}
A finite atlas and a subordinate partition of unity reduce
\eqref{eq:moser}--\eqref{eq:katoponce} to $\R^2$, where they are the
Moser product rule and the Kato--Ponce commutator estimate at real
order \cite[Chapter~13]{taylor2011}; the localization only adds
commutators of lower order, which interpolation absorbs.  Finally, the proof of Lemma~\ref{lem:lineartame} involves the order
$s$ only through \eqref{eq:moser}, \eqref{eq:katoponce}, and the
interpolation inequality
\(
 \|u\|_{H^{s+1}}
 \le\|u\|_{H^{s+2}}^{\,s/(s+1)}\,\|u\|_{H^{1}}^{\,1/(s+1)},
\)
all valid for real $s\ge2$, hence the same proof gives the statement with $p$ in place of $s$.
\end{proof}

\begin{proof}[Proof of the path-regularity bound \eqref{eq:pathregularity}]
Step 1: the bound for $h$.
By Proposition~\ref{prop:pathreg} with $(k,\alpha)=(4,\alpha_0)$,
each $h(t)$ lies in $C^{6,\alpha_0}(\sph)\subset C^4(\sph)$ and
satisfies \eqref{eq:uniformbounds}.  Moreover, $h(t)$ solves
\eqref{eq:perturbeddata} with $\psi=\rho_t$, where $\rho_t\ge4\delta$
and $\|\rho_t\|_{C^{2,\alpha}}\le\Lambda_*$ by
\eqref{eq:deltalambda}--\eqref{eq:alphachoice}.  Hence
Lemma~\ref{lem:apriori} applies with
$(\psi_-,\psi_+)=(\delta,\Lambda_*)$ and gives
$\|h(t)\|_{C^{4,\alpha}}\le M_1$ by \eqref{eq:corridorA}, as well as,
by \eqref{eq:tame},
\begin{align*}
 \sup_{t\in[0,1]}\|h(t)\|_{H^{s+2}}
 \le C\Bigl(1+\sup_{t\in[0,1]}\|\rho_t\|_{H^{s}}\Bigr)
 \le C\bigl(1+\|\rho\|_{H^{s}}\bigr).
\end{align*}

Step 2: the bound for $v$.
By \eqref{eq:velocity}, $v(t)\perp\mathcal N$ solves
$\mathcal L_{h(t)}v(t)=\rho-1$, where $\rho-1\in H^s(\sph)$ is
orthogonal to $\mathcal N$ by Lemma~\ref{lem:path}.  We apply
Lemma~\ref{lem:lineartame} with $h=h(t)$ and $f=\rho-1$; its
hypotheses hold uniformly in $t$, since $\kappa_0\sigma\le W[h(t)]\le\Lambda_0\sigma$ by \eqref{eq:uniformbounds},
$\|h(t)\|_{C^{4,\alpha}}\le M_1$ by Step~1, and the inf-sup constant is
bounded below by Proposition~\ref{prop:infsup}(i).  The
second bound in \eqref{eq:lineartame} of Lemma~\ref{lem:lineartame}
gives
\begin{align*}
 \|v(t)\|_{H^{s+2}}
 \le C\bigl(\|\rho-1\|_{H^s}
 +\bigl(1+\|h(t)\|_{H^{s+2}}\bigr)\|\rho-1\|_{C^\alpha}\bigr)\le C
\end{align*}
uniformly in $t$ by Step~1.
\end{proof}

\begin{proof}[Proof of Lemma~\ref{lem:perturbation}]
The first identity in \eqref{eq:bsym} holds because $b_h$ is
linear in $W[h]$, hence in $h$.  For smooth $e,\varphi$,
\eqref{eq:greens} and $\mathcal L_e\varphi=\mathcal L_\varphi e$, which holds by
\eqref{eq:Lop} and the symmetry of $D$, give
\begin{align*}
b_e(\varphi,q)=-(\mathcal L_e\varphi,q)
=-(\mathcal L_\varphi e,q)=b_\varphi(e,q),
\end{align*}
and both sides are continuous on $C^2\times C^2\times H^1$, so
density gives the second identity.  Finally,
\begin{align*}
|b_\varphi(e,q)|
\le C\|\varphi\|_{C^2}\bigl(\|\nabla e\|_{L^2}\|\nabla q\|_{L^2}
+\|e\|_{L^2}\|q\|_{L^2}\bigr)
\le C\|\varphi\|_{C^2}\|e\|_{H^1}\|q\|_{H^1}. \tag*{\qedhere}
\end{align*}
\end{proof}

\begin{proof}[Proof of Lemma~\ref{lem:lineartame}]
Step 1: existence and the H\"older estimate.
By the inf-sup bound \eqref{eq:infsup} and the symmetry of $b_h$,
there is a unique $w\in\Hperp$ with $b_h(w,q)=-(f,q)$ for all
$q\in\Hperp$, and $\|w\|_{H^1}\le2\beta^{-1}\|f\|_{(\Hperp)'}$.
As $f\perp\mathcal N$ and $b_h(w,\nu_k)=0$, the identity extends to
$q=\nu_k$, so $\mathcal L_hw=f$ weakly.  The coefficients of
$\mathcal L_h$ are bounded in $C^{2,\alpha}$ by $C(M_0)$, so elliptic
$H^2$ regularity, the embedding $H^2(\sph)\hookrightarrow C^0(\sph)$,
and the Schauder estimate on a fixed finite atlas
\cite[Chapters~6 and~8]{gilbargtrudinger2001} give
\begin{align*}
 \|w\|_{C^0}&\le C\|w\|_{H^2}
 \le C\bigl(\|f\|_{L^2}+\|w\|_{H^1}\bigr)\le C\|f\|_{L^2},\\
 \|w\|_{C^{2,\alpha}}&\le C\bigl(\|f\|_{C^\alpha}+\|w\|_{C^0}\bigr)
 \le C\|f\|_{C^\alpha}.
\end{align*}

Step 2: the Sobolev estimate.
In a chart of the fixed atlas, $\mathcal L_hw=a^{ij}\partial_i\partial_jw
+b^i\partial_iw+cw$, where $a^{ij}$, $b^i=-a^{jk}\Gamma^i_{jk}$,
and $c$ are linear in $h$ and its derivatives up to order two (through
$C[h]$); in particular
$\|a^{ij}\|_{H^s}+\|b^i\|_{H^s}+\|c\|_{H^s}\le C\|h\|_{H^{s+2}}$.
Applying $J^s=(I-\Delta)^{s/2}$ to the localized equation gives
\begin{align*}
 \mathcal L_h(J^sw)=J^sf-[J^s,a^{ij}]\,\partial_i\partial_jw
 -[J^s,c]\,w-R,
\end{align*}
with $R$ the remaining commutators, with $b^i\partial_i$ and
with the cutoffs, which by \eqref{eq:moser}--\eqref{eq:katoponce} satisfy
$\|R\|_{L^2}\le C\|w\|_{H^{s+1}}+C\|h\|_{H^{s+2}}\|w\|_{C^2}$.  By
\eqref{eq:katoponce},
\begin{align*}
 \|[J^s,a^{ij}]\,\partial_i\partial_jw\|_{L^2}
 +\|[J^s,c]\,w\|_{L^2}
 \le C\|w\|_{H^{s+1}}+C\|h\|_{H^{s+2}}\|w\|_{C^2}.
\end{align*}
The $H^2$ estimate of Step~1 applied to $J^sw$, together with these
commutator bounds, gives
$\|w\|_{H^{s+2}}\le C\bigl(\|f\|_{H^s}+\|w\|_{H^{s+1}}
+\|h\|_{H^{s+2}}\|w\|_{C^2}\bigr)$; interpolating
$\|w\|_{H^{s+1}}\le\varepsilon\|w\|_{H^{s+2}}+C_\varepsilon\|w\|_{H^1}$
and absorbing the first term yields
\begin{align}\label{eq:localtame}
 \|w\|_{H^{s+2}}
 \le C\Bigl(\|f\|_{H^s}+\|w\|_{H^1}
 +\bigl(1+\|h\|_{H^{s+2}}\bigr)\|w\|_{C^{2,\alpha}}\Bigr).
\end{align}
This a priori argument presumes $w\in H^{s+2}$; the standard smoothing
argument \cite[Chapter~14]{taylor2011} removes this assumption and
shows that $w$ lies in $H^{s+2}$.
\end{proof}

\begin{proof}[Proof of Lemma~\ref{lem:discreteinfsup}]
Step 1: exact coefficients.
Since $b_{h(t)}$ is only coercive up to a lower-order term, the
inf-sup bound of Proposition~\ref{prop:infsup} does not pass to the
subspace $V_L^\perp$ directly; we follow \cite{schatz1974}.
Fix $t$, write $b=b_{h(t)}$, and for $v\in V_L^\perp$ set
$\vartheta:=\sup_{q\in V_L^\perp\setminus\{0\}}|b(v,q)|/\|q\|_{H^1}$.
Lemma~\ref{lem:garding} gives
$\kappa_0\|v\|_{H^1}^2\le\vartheta\|v\|_{H^1}+C\|v\|_{L^2}^2$.
To control $\|v\|_{L^2}$ by $\vartheta$, let $z\in\Hperp$ solve
$b(z,q)=(v,q)$ for all $q\in\Hperp$, so that $\|z\|_{H^2}\le C\|v\|_{L^2}$ uniformly in $t$
by Proposition~\ref{prop:infsup} and elliptic regularity.  Since
$P_Lz\in V_L^\perp$, the symmetry of $b$ and \eqref{eq:approx} give
\begin{align*}
 \|v\|_{L^2}^2
 &=b(v,z-P_Lz)+b(v,P_Lz)\\
 &\le C\|v\|_{H^1}\|z-P_Lz\|_{H^1}+\vartheta\|P_Lz\|_{H^1}\\
 &\le C\bigl(L^{-1}\|v\|_{H^1}+\vartheta\bigr)\|v\|_{L^2},
\end{align*}
hence
$\kappa_0\|v\|_{H^1}^2\le\vartheta\|v\|_{H^1}+CL^{-2}\|v\|_{H^1}^2+C\vartheta^2$.
For $L\ge L_0$, absorption and Young's inequality yield
$\vartheta\ge2\beta_0\|v\|_{H^1}$ with $\beta_0>0$ independent of $t$.

Step 2: perturbed coefficients.
By the form bound in Step~4 of the proof of
Proposition~\ref{prop:infsup}, $b_h$ and $b_{h(t)}$ differ by at most
$2\varepsilon_0\|v\|_{H^1}\|q\|_{H^1}$, so $\varepsilon_0\le\beta_0/2$
gives \eqref{eq:discreteinfsup}.

Step 3: the load and the constraint.
With $\ell_L(q):=Q_L(fq)$, testing with $\nu_k$ fixes
$\lambda_k=-\tfrac{3}{4\pi}\ell_L(\nu_k)$ as for
\eqref{eq:lambdaformula}, and on $V_L^\perp$ the system reduces to
$b_h(v,q)=-\ell_L(q)$, uniquely solvable by \eqref{eq:discreteinfsup}.
As $V_L=V_L^\perp\oplus\mathcal N$ and $b_h(v,\nu_k)=0$, this $v$
solves the full saddle system, and
\begin{align*}
 \|v\|_{H^1}+|\lambda|
 \le C\sup_{q\in V_L\setminus\{0\}}\frac{|\ell_L(q)|}{\|q\|_{H^1}}
 \le C\,Q_L(f^2)^{1/2}
\end{align*}
by the Cauchy--Schwarz inequality for the positive weights of $Q_L$
and $Q_L(q^2)=\|q\|_{L^2}^2$.  With exact data integrals,
$\ell_L(q)=(f,q)$ and the same argument gives $\|f\|_{L^2}$.
\end{proof}

The Ricci identity, in the convention of the proof of
Lemma~\ref{lem:divfree}, also yields the Simons identity used in
Step~1 of the proof of Lemma~\ref{lem:apriori}: for $h\in C^4(\sph)$,
$A=W[h]$, and $H=\tr A$,
\begin{align}\label{eq:simons}
 \Delta A=\nabla^2H+2A-H\sigma .
\end{align}
The Codazzi property of $A$ and the divergence identity
$\nabla^kA_{kj}=\nabla_jH$ of Lemma~\ref{lem:divfree} give
\begin{align*}
 \Delta A_{ij}
 =\nabla^k\nabla_kA_{ij}
 =\nabla^k\nabla_iA_{kj}
 =\nabla_i\nabla_jH+[\nabla^k,\nabla_i]A_{kj}.
\end{align*}
The commutator acts on each slot of $A$ through
$R_{pik}{}^{m}=\sigma_{pk}\delta_i^m-\delta_p^m\sigma_{ik}$, and
since $\sigma^{kp}\sigma_{pk}=2$,
\begin{align*}
 \sigma^{kp}R_{pik}{}^{m}A_{mj}=2A_{ij}-A_{ij}=A_{ij},
 \qquad
 \sigma^{kp}R_{pij}{}^{m}A_{km}=A_{ij}-H\sigma_{ij}.
\end{align*}
Their sum is $[\nabla^k,\nabla_i]A_{kj}=2A_{ij}-H\sigma_{ij}$, which
proves \eqref{eq:simons}.

\begin{proof}[Proof of the estimate \eqref{eq:tame} in Lemma~\ref{lem:apriori}]
Let $X\in\{X_1,X_2,X_3\}$.  By \eqref{eq:killing}, $w=Xh$ solves
$\mathcal L_hw=X\psi$.  The a priori estimate \eqref{eq:localtame} applies to $w$, since its
derivation uses only the equation $\mathcal L_hw=f$ and the bound
$\|h\|_{C^{4,\alpha}}\le M_1$ of \eqref{eq:corridorA}, not the
inf-sup hypothesis of Lemma~\ref{lem:lineartame}.

Let $4\le p'\le p+1$ and $h\in H^{p'}(\sph)$.  Then $w\in H^{p'-1}$
and $X\psi\in H^{p'-2}$, as $p'-1\le p$, so the smoothing argument
cited after \eqref{eq:localtame} gives $w\in H^{p'}$.  Since
$\|w\|_{H^1}+\|w\|_{C^{2,\alpha}}\le CM_1$ by \eqref{eq:C4alphabound},
\eqref{eq:localtame} at order $p'-2$ yields
$\|Xh\|_{H^{p'}}\le C\bigl(1+\|\psi\|_{H^p}+\|h\|_{H^{p'}}\bigr)$.
As the $X_j$ are skew-adjoint, commute with $\Delta$, and satisfy
$\sum_jX_j^2=\Delta$, \eqref{eq:sobolev} gives
$\|h\|_{H^{p'+1}}^2=\|h\|_{H^{p'}}^2+\sum_j\|X_jh\|_{H^{p'}}^2$, hence
\begin{align*}
 \|h\|_{H^{p'+1}}\le C\bigl(1+\|\psi\|_{H^p}+\|h\|_{H^{p'}}\bigr).
\end{align*}
Starting from $\|h\|_{H^4}\le C\|h\|_{C^{4,\alpha}}\le CM_1$, we apply
this with $p'=4,5,\dots,\lceil p\rceil$ and finally with $p'=p+1$.
These at most $\lceil p\rceil-2$ steps give \eqref{eq:tame}, with a
constant depending only on $\psi_-,\psi_+,\alpha,p$.
\end{proof}

\subsection{Verification of the umbilic phase}\label{app:phase}

Let $a$ be the angle parametrizing the circle $\Gamma$ of Section~\ref{subsec:expU}, $r=0.08$
its radius, and $\theta_0\approx0.8355$ the polar angle of its
center.  In the frame tangent and normal to $\Gamma$ within $\sph$,
let $\widetilde q$ be defined as in Section~\ref{subsec:expU} with this
frame in place of $(e_\theta,e_\varphi)$; then $|\widetilde q|=|q|$.
The Cartesian components of $\Gamma$ and of this frame are
trigonometric polynomials of degree one in $a$, and those of
$W[h_L]$ have degree $L+2$ by Step~2 of the proof of
Lemma~\ref{lem:Winfty}, so $\widetilde q$ has degree at most
$L+4=28$.

At $a_j=2\pi j/720$, $j=0,\ldots,719$, the $L=24$ coefficients give
$\min_j|q(a_j)|>0.067$ and $\max_j|q(a_j)|<0.080$.  Since every
$a$ lies within $\pi/720$ of a sample, Bernstein's inequality yields
\begin{align*}
 \|\widetilde q\|_\infty&\le\frac{0.080}{1-28\pi/720}<0.092,\\
 \|\partial_a\widetilde q\|_\infty&\le28\cdot0.092<2.6,\\
 \min_\Gamma|q|&\ge0.067-2.6\,\pi/720>0.055.
\end{align*}
Let $\beta$ be the angle from the spherical frame to the circle
frame.  The geodesic curvature of $\Gamma$ and the connection form
of the spherical frame give
$|\beta'|\le\cos r+\sin r\cot(\theta_0-r)<1.1$, and a frame rotation
by $\beta$ shifts $\arg q$ by $2\beta$.  Hence
\begin{align*}
 |\partial_a\arg q|
 \le\frac{\|\partial_a\widetilde q\|_\infty+2\,|\beta'|\,\|\widetilde q\|_\infty}{\min_\Gamma|q|}
 <51,
  \quad
 \bigl|\arg q(a_{j+1})-\arg q(a_j)\bigr|<51\cdot\frac{2\pi}{720}
 <0.45<\pi ,
\end{align*}
between consecutive samples.  Thus $\Gamma$ is zero-free, the
unwrapping is unambiguous, and the sampled winding, which equals one,
is the true winding.

\bigskip
{\scshape Department of Applied Mathematics, The Hong Kong Polytechnic University, Hung Hom, Kowloon, Hong Kong}\par
\emph{Email address}: \href{mailto:buyang.li@polyu.edu.hk}{\texttt{buyang.li@polyu.edu.hk}}\par
\medskip
{\scshape Department of Mathematics and Department of Physics, Duke University, Durham, NC 27707, USA}\par
\emph{Email address}: \href{mailto:jian-guo.liu@duke.edu}{\texttt{jian-guo.liu@duke.edu}}\par
\medskip
{\scshape Department of Applied Mathematics, The Hong Kong Polytechnic University, Hung Hom, Kowloon, Hong Kong}\par
\emph{Email address}: \href{mailto:claire.tang@polyu.edu.hk}{\texttt{claire.tang@polyu.edu.hk}}\par

\end{document}